\documentclass{amsart}
\usepackage{amssymb}

\usepackage{graphicx}
\usepackage{amscd}

\theoremstyle{plain}

\newtheorem{definition}{Definition}

\newtheorem{lemma}{Lemma}

\newtheorem{proposition}{Proposition}
\newtheorem{remark}{Remark}

\newtheorem{theorem}{Theorem}
\numberwithin{equation}{section}

\begin{document}
\title[weakly stationary-harmonic multiple-valued functions]
{Interior continuity of two-dimensional weakly stationary-harmonic
multiple-valued functions}
\author{Chun-Chi Lin}
\address{Department of Mathematics, National Taiwan Normal University}
\email{chunlin@math.ntnu.edu.tw, chunlin@mis.mpg.de}
%\thanks{}
\date{May 09, 2013}

\subjclass[2000]{Primary 49Q20; Secondary 28A75, 35A15}
\keywords{Almgren's big regularity paper, multiple-valued functions,
weakly stationary-harmonic, interior continuity}

\begin{abstract}
In his big regularity paper, Almgren has proven the regularity theorem for mass-minimizing integral 
currents.
One key step in his paper is to derive the regularity of Dirichlet-minimizing
$\mathbf{Q}_{Q}(\mathbb{R}^{n})$-valued functions in the Sobolev space
$\mathcal{Y}_{2}(\Omega, \mathbf{Q}_{Q} (\mathbb{R}^{n}))$,
where the domain $\Omega$ is open in $\mathbb{R}^{m}$.
In this article, we introduce the class of weakly stationary-harmonic
$\mathbf{Q}_{Q} (\mathbb{R}^n)$-valued functions.
These functions are the critical points of Dirichlet integral under
smooth domain-variations and range-variations.
We prove that if $\Omega$ is a two-dimensional domain in $\mathbb{R}^{2}$
and 
$f\in\mathcal{Y}_{2}\left(\Omega,\mathbf{Q}_{Q}(\mathbb{R}^{n})\right)$
is weakly stationary-harmonic,
then $f$ is continuous in the interior of the domain $\Omega$.

\end{abstract}

\maketitle

\baselineskip=15 pt plus 1pt minus .5pt

\section{Introduction}

In his big regularity paper \cite{Alm00},
Almgren proved the regularity theorem for mass-minimizing integral currents.
More precisely, Almgren showed that any mass-minimizing integral current is smooth except on a singular subset of co-dimension two.
One key step in \cite{Alm00} is to derive the regularity of Dirichlet-minimizing
multiple-valued functions in the Sobolev space
$\mathcal{Y}_{2}(\Omega, \mathbf{Q}_{Q} (\mathbb{R}^{n}))$
of multiple-valued functions (see (\ref{eq:Def-Y_2}) for its definition).
This is a key step because Almgren
used this class of multiple-valued functions
to approximate mass-minimizing
integral currents, whose regularity therefore inherits that of Dirichlet-minimizing
multiple-valued functions.
In \cite{Alm00}, any Dirichlet-minimizing multiple-valued function,
$f:\Omega\subset\mathbb{R}^{m}\rightarrow\mathbf{Q}_{Q}(\mathbb{R}^n)$,
is shown to be H\"{o}lder continuous at any interior point of $\Omega$
and smooth outside a closed subset $\Sigma$, whose co-dimension is least two
(in the sense of Hausdorff measure).
Moreover, as $\Omega\subset\mathbb{R}^{2}$,
the closed subset $\Sigma$ consists of isolated points.
Dirichlet-minimizing multiple-valued functions were further investigated 
in \cite{Zhu08} by Zhu as $\Omega\subset\mathbb{R}^2$ 
and 
in \cite{Spadaro10} by Spadaro as $\Omega\subset\mathbb{C}^m$. 
The reader is also referred to a recent article by
De Lellis and Spadaro in \cite{DS11},
which makes Almgren's theory of multiple-valued functions
more accessible and provides some new points of view on the subject,
e.g., intrinsic theory of the metric space $\mathbf{Q}_{Q}(\mathbb{R}^n)$.

Although the theory of multiple-valued functions was originally developed for
the purpose of studying the regularity of mass-minimizing integral currents
in \cite{Alm00},
we found that the theory itself is interesting enough
and deserves further investigation
as an independent topic in calculus of variations.
In \cite{Matt83}, Mattila investigated a class of elliptic variational integrals
of multiple-valued functions
(interior H\"{o}lder continuity was derived in $2$-dimensional case).
Along the direction of \cite{Matt83},
De Lellis, Focardi and Spadaro in \cite{DFS11} further
characterized the semicontinuity of certain elliptic integrals of
multiple-valued functions.
In \cite{GZ08}, Goblet and Zhu studied the regularity of
Dirichlet nearly minimizing multiple-valued functions.
On the other hand, Almgren's multiple-valued functions 
have also been used to formulate certain (stable) branched minimal surfaces 
and minimal hypersurfaces by Rosales, Simon, Wickramasekera, et al 
(e.g., see \cite{Rosales10}, \cite{Rosales11}, \cite{SW07}, \cite{Wick08}). 
Their research projects provide an approach to investigate 
more details of the local behavior around singularities of minimal surfaces 
and minimal hypersurfaces, 
which are not area-minimizing but could be formulated as 
two-valued minimal graphs.

In this article, we propose studying a bigger class of multiple-valued functions,
called \texttt{weakly stationary-harmonic} multiple-valued functions,
which are the critical points of Dirichlet integral
with respect to smooth domain-variations and range-variations
in the Sobolev space
$\mathcal{Y}_{2}(\Omega, \mathbf{Q}_{Q} (\mathbb{R}^n))$. 
At the first glance, the reader might consider this variational problem 
as an extension of theory for stationary harmonic maps. 
Indeed, this is what we thought at the beginning as we started this project.
However, there are several difficulties when one tries to apply the methods in the (partial) 
regularity theory for stationary harmonic maps.
The main difficulty comes from the fact that the metric space 
$\mathbf{Q}_{Q}(\mathbb{R}^n)$
could only be embedded into Euclidean spaces
as a bi-Lipschitz (polyhedral cone) submanifold,
which is not negatively curved and lacks sufficient smoothness.
Thus, we could not follow the typical definitions of stationary harmonic maps
and the approach studying their (partial) regularity in literature here.
We prove that if
$f\in\mathcal{Y}_{2}(\Omega, \mathbf{Q}_{Q} (\mathbb{R}^n))$
is weakly stationary-harmonic and $\Omega\subset\mathbb{R}^{2}$,
then $f$ is continuous in the interior of the domain $\Omega$.

There are two parts in the proof of our main result in Theorem \ref{thm:main}.
In the first part, we apply the domain-variations of $f$
to show that the Hopf differential of
$\boldsymbol\xi_{0}\circ f$
(denoted as $\Phi:\Omega\rightarrow\mathbb{C}$)
is a holomorphic function,
where
$\boldsymbol\xi_{0}:\mathbf{Q}_{Q} (\mathbb{R}^{n})
\rightarrow\mathbb{R}^{nQ}$
is any Lipschitz map defined in (\ref{eq:xi_0-Almgren}).
This trick has been used in the theory of $2$-dimensional harmonic maps 
with very general target spaces,
e.g., non-positively curved metric spaces (cf. \cite{Schoen84}).
Then, by applying another trick in Gr\"{u}ter's paper \cite{Gr86},
there is an induced harmonic function
$h:\Omega \rightarrow \mathbb{R}^{2}$
so that
$(\boldsymbol\xi_{0} \circ f, h):\Omega\rightarrow \mathbb{R}^{nQ+2}$
is weakly conformal.
This weakly conformal condition puts us in a position to adapt the method in
Gr\"{u}ter's paper \cite{Gr81}
(studying the regularity of weak conformal H-surfaces)
for the class of multiple-valued functions of this article in the second part.

In the second part, one key step is to establish a ``global" monotonicity formula
(see Step 2 of the proof of Lemma \ref{lem:prethm}) from
proper range-variations of $f$.
This part contains the main difficulty of this article
as we apply Gr\"{u}ter's approach in \cite{Gr81} to our case.
Notice that,
one could not simply perturb a multiple-valued function
by adding a test function belonging to a proper class of functions
(just like in PDE theory)
or by further projections into a smooth submanifold
(just like in the theory of stationary harmonic maps).
This difficulty is due to the lack of structure for algebraic operations
(e.g., addition and subtraction) in $\mathbf{Q}_{Q} (\mathbb{R}^n)$ as $n\ge 2$.
Even in the case of $n=1$,
where a natural subtraction between two members in
$\mathbf{Q}_{Q} (\mathbb{R})$ does exist (see the Preliminaries),
one still needs to be careful in perturbations of multiple-valued functions
so that Definition \ref{def:WH} is fulfilled.
A simple example is to naively define a family of multiple-valued functions from
subtraction between a
$\mathbf{Q}_{2} (\mathbb{R})$-valued function $f(x):=[\![x]\!]+[\![-x]\!]$,
where $x\in(-1,1)$, and the member $[\![-1]\!]+[\![1]\!]$
by $f^{t}(x):=[\![(1-t)x+t(1-x)]\!]+[\![-(1-t)x+t(1+x)]\!]$.
As $t>0$, it is easy to verify that there is no point with multiplicity $2$ remained
in the value of the continuous function $f^{t}$.
According to Definition \ref{def:WH},
one only allows (parametrized) diffeomorphisms
(or bi-Lipschitz homeomorphisms)
$\psi_{x}:=\psi(x,\cdot):\mathbb{R}^n\rightarrow\mathbb{R}^n$,
where
$\psi\in C^{\infty}_{c}(\Omega\times\mathbb{R}^n,\mathbb{R}^n)$,
in the range-variations.
Thus,
the point with multiplicity $2$ of the continuous
$\mathbf{Q}_{2} (\mathbb{R})$-valued function $f$ 
still remains in the class of range-variations given in Definition \ref{def:WH}.
Therefore, the family of perturbed
$\mathbf{Q}_{2}(\mathbb{R})$-valued functions $f^t$ can't be generated from
the definition of range-variations.

To overcome the difficulty in range-variations,
we give a method to build up the so-called \emph{ nested admissible} closed balls
(see Definition \ref{def:Nested_Adm_Closed_Ball})
for arbitrarily chosen member in $\mathbf{Q}_{Q} (\mathbb{R}^n)$.
Namely, for any member $q\in\mathbf{Q}_{Q} (\mathbb{R}^n)$,
there exists a sequence $q=q^{(0)}, q^{(1)}, ..., q^{(L)}=Q[\![ a ]\!]$ for some
$L\in\mathbb{Z}_{+}$ and $a\in\mathbb{R}^n$ with strictly decreasing
$\text{card}\left(\text{spt}\left(q^{(i)}\right)\right)$ as $i$ increases
such that
\begin{equation*}
\cdots\subset\subset
\mathbb{B}^{\mathbf{Q}}_{\rho_{i}}(q^{(i)})
\subset\subset\mathbb{B}^{\mathbf{Q}}_{\sigma_{i}}(q^{(i)})
\subset\subset\mathbb{B}^{\mathbf{Q}}_{\rho_{i+1}}(q^{(i+1)})
\subset\subset\cdots
\end{equation*}
for some sequence
\begin{equation*}
0=\rho_{0}<\sigma_{0}<\rho_{1}<\sigma_{1}<\cdots<\rho_{L}<\sigma_{L}:=\infty
\text{.}
\end{equation*}
One observes that there is a natural subtraction between this arbitrarily chosen member
$q$ and any member in its \emph{admissible} closed balls (or neighborhoods).
Then we show that the monotonicity formula can be established by
the perturbations associated with this type of \emph{admissible} balls
on $\mathbb{B}^{\mathbf{Q}}_{\rho}(q^{(i)})$
(see Definition \ref{def:Adm_Closed_Ball}), where
$\rho$ is roughly between $\rho_{i}$ and $\sigma_{i}$.
The monotonicity formula would be used to bound
$\sigma_{i}-\rho_{i}$ from above by Dirichlet integral of
$(\boldsymbol\xi_{0}\circ f, h)$ for each $i$.
Note that $\rho_{i+1}-\sigma_{i}$ is bounded from above by a constant
depending only on $n$ and $Q$
(see Proposition \ref{prop:Nested_Adm_Closed_Ball}).
This allows us to establish a ``global" monotonicity formula
in the proof of our key lemma (see Step 3 in the proof of Lemma \ref{lem:prethm}). 
Here, ``global" means that the formula is not restricted to the 
radius of an admissible closed ball that one derives the formula
but is extensible to a bigger admissible closed ball
(with changes of some constants depending only on $n$ and $Q$).

Finally, 
in the end of this article, 
we show the proof of interior continuity for weakly stationary-harmonic 
$\mathbf{Q}_{Q} (\mathbb{R}^{n})$-valued functions, 
using a contradiction argument 
by applying the key lemma (i.e., Lemma \ref{lem:prethm}) 
and Courant-Lebesgue Lemma (see Lemma \ref{lem:CL}).

\section{Preliminaries}
In this section, we collect some notations and results from
\cite{Alm00}, \cite{FH1} and \cite{Gr81} to keep this article short
and self-contained. 

Denote by $[\![p_i]\!]$ the Dirac measure at the point $p_i\in \mathbb{R}^{n}$.  
For a given positive integer $Q$, 
let
$\mathbf{Q}_{Q}(\mathbb{R}^n)
:=\{\sum\limits_{i=1}^{Q}[\![p_{i}]\!]:p_{i}\in \mathbb{R}^n \}$,
where $p_{i}$, $p_{j}$ are not necessarily distinct as $i\neq j$.
For our convenience, as we use the notation
$p=\sum\limits_{j=1}^{J} \ell_j [\![p_j]\!]$ for any member in
$\mathbf{Q}_{Q}(\mathbb{R}^{n})$, 
$p_1,...,p_J$ are distinct points in
$\mathbb{R}^n$ 
and $\ell_j$ is the multiplicity of $p_j$ for each $j$ 
(therefore $\sum\limits_{j=1}^{J} \ell_j=Q$). 
In \cite{Alm00},
$\mathbf{Q}_{Q}(\mathbb{R}^n)$
is shown to be a metric space by associated with proper distance functions.
One is the so-called flat metric,
see {\cite[Chap.1]{Alm00}}.
The other one, 
corresponding to the standard Euclidean distance function, 
is given by
\begin{equation*}
\mathcal{G}(p,q)
:= {\inf}
\left\{
\left(
\overset{Q}{\underset{i=1}{\sum }}
\left| p_{i}-q_{\sigma (i)}\right| ^{2}
\right)^{1/2}:
\text{$\sigma$ is a permutation of $\{1,..., Q\}$}
\right\}
\end{equation*}
where
$p=\underset{i=1}{\overset{Q}{\sum }}[\![p_{i}]\!]$
and
$q=\underset{j=1}{\overset{Q}{\sum }}[\![q_{j}]\!]$. 
Let $\Omega\subset\mathbb{R}^{m}$ be an open set. 
As $Q\ge 2$, 
$f:\Omega\subset\mathbb{R}^{m}\rightarrow \mathbf{Q}_{Q}(\mathbb{R}^n)$
is called a \texttt{multiple-valued function}
or more precisely a 
$\mathbf{Q}_{Q}(\mathbb{R}^n)$-\texttt{valued function},
which is usually denoted by 
$f(x)=\underset{i=1}{\overset{Q}{\sum }}[\![f_{i}(x)]\!]$. 
Therefore, the support of $f(x)$, $\text{spt}\left(f(x)\right)$,
consists of $Q$ unordered points in
$\mathbb{R}^{n}$ for all $x\in \Omega$.
In fact, Almgren showed in {\cite[Section 1.1]{Alm00}}
and {\cite[Section 1.2]{Alm00}} that the metric space
$\mathbf{Q}_{Q}(\mathbb{R}^n)$ is
bi-Lipschitz correspondence with a $nQ$-dimensional polyhedral cone
$\mathbf{Q}^{\ast}$ in a Euclidean space
$\mathbb{R}^{P(n,Q) \cdot Q}$.
More precisely, there exists a positive integer $P(n,Q)>n$ such that
for each fixed $\alpha=1,..., P(n,Q)$,
there corresponds a straight line $L_{\alpha}$ in $\mathbb{R}^n$
and an orthogonal projection, denoted as $\Pi_{\alpha}\in O^{\ast }(n,1)$,
into the $\alpha$-th straight line in $\mathbb{R}^{n}$.
For each $\alpha\in\{1,2,...,n\}$, 
$L_{\alpha}$ is chosen to be
the coordinate axis of $\mathbb{R}^n$.
Therefore, for any $y=(y^1,...,y^n)\in\mathbb{R}^n$
and $\alpha \in\{1,..., n\}$,
we have $\Pi_{\alpha}(y)=y^\alpha \in\mathbb{R}$.
Thus, for each $\alpha\in\{1,...,n\}$,
$\Pi_{\alpha}$ induces the map
$(\Pi_{\alpha})_{\#}(\cdot):\mathbf{Q}_{Q}(\mathbb{R}^{n})
\rightarrow\mathbf{Q}_{Q}(\mathbb{R})$
defined by
\begin{equation}
(\Pi_{\alpha})_{\#}\left(\underset{j=1}{\overset{Q}{\sum}}[\![q_{j}]\!]\right)=
\underset{j=1}{\overset{Q}{\sum }}[\![q_{j}^{\alpha}]\!] 
\label{eq:(Pi_(alpha))_(sharp)}
\end{equation}
where
$q_{j}^{\alpha}$ is the $\alpha$-th component of
$q_{j}\in\mathbb{R}^n$,
and
\begin{equation}
\xi_{\alpha}(\cdot):=
{\xi}(\Pi_{\alpha},\cdot ): \mathbf{Q}_{Q} (\mathbb{R}^{n})
\rightarrow
\left\{
(s_{1}, s_{2}, \cdot \cdot \cdot, s_{Q}):
s_{1}\le s_{2}\le\cdot\cdot\cdot\le s_{Q}
\right\}
\subset \mathbb{R}^{Q}
\text{.}
\label{eq:Xi_0}
\end{equation}
Therefore,
\begin{equation*}
\xi_{\alpha}(q)=
(q_{\sigma(1)}^{\alpha}, \cdot \cdot \cdot, q_{\sigma(Q)}^{\alpha})
\text{,}
%\label{eq:Xi_alpha}
\end{equation*}
for some $\sigma\in \mathcal{P}_{Q}$: the permutation group of
$\{1,...,Q\}$.
Note that, from (\ref{eq:Xi_0}), it is easy to verify
\begin{equation}
\left| \xi_{\alpha}(p)-\xi_{\alpha}(q) \right|
=\left|{\xi}(\Pi_{\alpha},p) - {\xi}(\Pi_{\alpha},q) \right|
= \mathcal{G} ( (\Pi_{\alpha})_{\#}(p) , (\Pi_{\alpha})_{\#}(q) )
\text{.}
\label{eq:Q(R)-topology}
\end{equation}
For each fixed coordinates of $\mathbb{R}^n$,
Almgren introduce the Lipschitz correspondence 
\begin{equation}
\boldsymbol\xi_{0}(\cdot)
:=\left(\xi(\Pi_{1}, \cdot)\cdot\cdot\cdot\xi(\Pi_{n}, \cdot)\right)
:\mathbf{Q}_{Q}(\mathbb{R}^{n})\rightarrow \mathbb{R}^{n Q}
\label{eq:xi_0-Almgren}
\end{equation}
with $Lip(\boldsymbol\xi _{0})=1$.
The reader may verify from (\ref{eq:Q(R)-topology}) that
\begin{equation*}
|\boldsymbol\xi_{0}(p)-\boldsymbol\xi_{0}(q)|^2
= \sum\limits_{\alpha=1}^{n} |\xi_{\alpha}(p)-\xi_{\alpha}(q)|^2
\le \mathcal{G}\left( p, q \right)^2
\text{.}
\end{equation*}
Moreover, $\boldsymbol\xi_{0}$ is not injective
unless $n=1$ (or see {\cite[Theorem 1.2]{Alm00}}).
However,
by introducing more distinct orthogonal projections into the straight lines
$L_{\alpha}$, $\alpha\in\{n+1,...,P(n,Q)\}$
(as defined in (\ref{eq:Xi_0})),
Almgren showed that
\begin{equation*}
\boldsymbol\xi(\cdot):= \xi(\Pi_{1},\cdot)
\Join \cdot \cdot \cdot \Join
\xi(\Pi_{P},\cdot):\mathbf{Q}_{Q}(\mathbb{R}^{n})
\rightarrow \mathbf{Q}^{\ast}
\subset \mathbb{R}^{N}
\end{equation*}
is then a bi-Lipschitz correspondence,
where
$\mathbf{Q}^{\ast}
:=\boldsymbol\xi\left(\mathbf{Q}_{Q}(\mathbb{R}^{n})\right)$
and $N=P(n,Q) Q$.
Furthermore,
both
$Lip(\boldsymbol\xi)$ and $Lip(\boldsymbol\xi^{-1})$
are positive numbers depending only on $n$ and $Q$.

Besides the bi-Lipschitz correspondence
$\boldsymbol\xi$ introduced in
\cite{Alm00},
there is a modified bi-Lipschitz and locally (or infinitesimally)
equidistant correspondence introduced by Brian White.
The reader could find it in literature from the article of Sheldon Chang
(see {\cite[p. 706]{Chang88}})
or from \cite{DS11} for more details.
The modified correspondence is constructed by choosing
$P=P(n,Q)$ distinct orthonormal coordinate bases 
(by rotating the orthonormal coordinates of $\mathbb{R}^n$),
by taking the orthogonal projections $\Pi_{1},...,\Pi_{P\cdot n}$
(as done in {\cite[Chapter 1.2]{Alm00}}) as a complete set of coordinate projections,
and by rescaling the resulting $\boldsymbol\xi$
under a proper scaling factor depending on $P(n,Q)$.

%One might need to revise the paragraph below, especially the notations of affine maps, linear maps, the L^2-norm of affine maps (or the derivatives),
%$\text{ap }D_{x} f$ or $\text{ap }D_{x_{0}} f$... etc.
%!!!!!!!!!!!!!!!!!!!!!!!!!!!!!!!!!!!!!!!!!!!!!!!!!!!!!!!!!!!!!!!!!!!!!!!!!
Denote an affine map
$A: \mathbb{R}^{m}\rightarrow\mathbb{R}^n$
by
$A(x)=A(x_0)+L(x-x_0)$,
where
$L\in \text{Hom}(\mathbb{R}^{m},\mathbb{R}^n)$ is the linear part.
Let 
$\text{A}(m, n)$
denote the set of affine maps from $\mathbb{R}^{m}$ to $\mathbb{R}^n$.
As $A\in\text{A}(m,n)$,
we let
\begin{equation*}
|A|:=
\left(\sum\limits_{i=1}^{m} \text{ }
\left| D_i L\right|^2
\right)^{1/2}\in\mathbb{R} 
\end{equation*}
where
$D_{i} L=\frac{\partial L}{\partial x_{i}}$.
A multiple-valued function
$\mathcal{A}:\mathbb{R}^m \rightarrow \mathbf{Q}_{Q}(\mathbb{R}^n)$
is said to be \texttt{affine} if there are affine maps
$A_{1},..., A_{Q}\in\text{A}(m, n)$
such that
$\mathcal{A}:=\sum\limits_{i=1}^{Q}[\![A_{i}]\!]$.
Let
\begin{equation}
|\mathcal{A}|:=
\left(\sum_{i=1}^{Q} \text{ } \left| A_i\right|^2\right)^{1/2}
\text{.}
\label{eq:affine_norm}
\end{equation}
Assume $\Omega \subset \mathbb{R}^{m}$ is an open set and
$x_0\in\Omega$,
then
$f:\Omega\rightarrow \mathbf{Q}_{Q}(\mathbb{R}^n)$ is called
\texttt{approximately affinely approximatable} at $x_{0}\in\Omega$
if there exists an affine function
$\mathcal{A}:\mathbb{R}^{m}\rightarrow\mathbf{Q}_{Q}(\mathbb{R}^n)$
such that
\begin{equation*}
\underset{x\rightarrow x_0}{\text{ap } \lim}
\frac{\mathcal{G}(f(x), \mathcal{A}(x))}{\mid x-x_{0}\mid}=0
\text{.}
\end{equation*}
Such multiple-valued function
$\mathcal{A}$ is uniquely determined and denoted by
$\text{ap }Af(x_0)$.
Thus, as
$f:\Omega\rightarrow \mathbf{Q}_{Q}(\mathbb{R}^n)$
is approximately affinely approximatable at $x_{0}\in\Omega$,
we write
\begin{equation*}
\text{ap }Af(x_0)= \sum\limits_{i=1}^{Q}[\![A_{i}(x_0)]\!]
\text{.}
\end{equation*}
In {\cite[Theorem 1.4 (3)]{Alm00}}, Almgren proved that if
$\boldsymbol\xi\circ f:\Omega\rightarrow\mathbb{R}^{N}$ is approximately differentiable at $x_{0}\in\Omega$ for a multiple-valued function
$f(x)=\underset{i=1}{\overset{Q}{\sum }}[\![f_{i}(x)]\!]\in
\mathbf{Q}_{Q}(\mathbb{R}^{n})$,
then $f$ is
\texttt{strongly approximately affinely approximatable} at $x_0$.
In other words, 
if $\boldsymbol\xi\circ f$ is approximate differentiable at $x_0$,
then $f_{i}(x_0)=f_{j}(x_0)$ implies that
\begin{equation*}
\text{ap } D_{x} f_{i}(x_0) = \text{ap } D_{x} f_{j}(x_0)
\text{.}
\end{equation*}
Furthermore,
$\text{ap } Af(x_0)$ is uniquely determined by
$f(x_0)$ and $\text{ap } D\left(\boldsymbol\xi\circ f\right)(x_0)$
and
\begin{equation*}
|\text{ap }Af(x_0)|
=|\text{ap }D\left(\boldsymbol\xi_{0}\circ f\right)(x_0)| 
\text{.}
%\label{eq:apAf=ap D(xi-f)-1}
\end{equation*}

In \cite{Alm00}, Almgren also introduced the space,
$\mathcal{Y}_{2}(\Omega,\mathbf{Q}_{Q}(\mathbb{R}^{n}))$
(still called Sobolev spaces for our convenience),
for the class of multiple-valued functions
$f:\Omega\subset\mathbb{R}^m\rightarrow\mathbf{Q}_{Q}(\mathbb{R}^{n})$.
Recall that one usually uses 
$W^{1, 2}(\Omega, \mathbb{R}^{N})$
to denote the Sobolev space of
$\mathbb{R}^{N}$-valued functions defined on $\Omega$
with their first order distributional partial derivatives being
$\mathcal{L}^m$ square summable over $\Omega$. 
Denote by 
$\mathbb{U}_{r}^{m}(x)$ the open ball of radius $r$ 
with center $x$ in $\mathbb{R}^m$; 
and 
by $\mathbb{U}_{r}(x)$ as $m=2$ for convenience. 
A function
$f \in W^{1, 2}(\Omega, \mathbb{R}^{N})$
is said to be \emph{strictly defined} if
$f(x)=y$ as $x\in \Omega$, $y\in \mathbb{R}^{N}$, and
\begin{equation*}
\underset{r\rightarrow 0}{\lim} r^{-m} \text{ }
\int_{\mathbb{U}_{r}^{m}(x)}\left| f (z)-y \right|\text{ }d\mathcal{L}^m z
=0
\text{.}
\end{equation*}
In fact, any
$f \in W^{1, 2}(\Omega, \mathbb{R}^{N})$
agrees with a strictly defined
$g \in W^{1, 2}(\Omega, \mathbb{R}^{N})$
$\mathcal{L}^m $ a.e. on $\Omega$
(see {\cite[Appendix 1.2]{Alm00}} or {\cite[p. 592]{Matt83}}).
In \cite{Alm00}, the Sobolev space for multiple-valued functions
is defined by
\begin{equation}
\mathcal{Y}_{2}(\Omega ,\mathbf{Q}_{Q}(\mathbb{R}^{n}))
:=
\{
f:\Omega\subset\mathbb{R}^m\rightarrow\mathbf{Q}_{Q}(\mathbb{R}^{n})
: \boldsymbol\xi\circ f \in W^{1, 2}(\Omega, \mathbb{R}^{N})
\}
\text{.}
\label{eq:Def-Y_2}
\end{equation}
We say that $f$ is \emph{strictly defined} if
$\boldsymbol\xi\circ f$ is strictly defined.
Suppose
$f\in \mathcal{Y}_{2}(\Omega ,\mathbf{Q}_{Q}(\mathbb{R}^{n}))$,
then by {\cite[Theorem 2.2]{Alm00}}, 
\begin{equation}
\left| \text{ap }D \left(\boldsymbol\xi_{0}\circ f(x)\right) \right|
=\left| \text{ap }Af(x)\right|
\label{eq:apAf=ap D(xi-f)-2}
\end{equation}
a.e. $x\in\Omega$. 
From {\cite[Definition 2.1]{Alm00}} and {\cite[Theorem 2.2]{Alm00}},
the Dirichlet integral of a multiple-valued function
$f\in \mathcal{Y}_{2}(\Omega, \mathbf{Q}_{Q}(\mathbb{R}^{n}))$
over an open set
$\Omega$ is given by
\begin{equation}
Dir(f; \Omega)
=
\underset{\Omega}{\int}
\left| \text{ap } A f(x) \right|^{2}
\text{ }
d\mathcal{L}^{m} x
=
\underset{\Omega}{\int}
\left| \text{ap } D \left( \boldsymbol\xi_{0}\circ f(x)\right) \right|^{2}
\text{ }
d\mathcal{L}^{m} x
\text{.}
\label{eq:Dir-Almgren}
\end{equation}
Note that,
since the norm of any affine map defined in (\ref{eq:affine_norm})
is independent of the choice of coordinates of $\mathbb{R}^n$,
we have
\begin{equation}
|\text{ap } D (\boldsymbol\xi_0\circ f )|
=
|\text{ap } D (\overset{\sim}{\boldsymbol\xi_0}\circ f )| 
\label{eq:eq_xi_0-f}
\end{equation}
for any two distinct Lipschitz correspondences
$\boldsymbol\xi_0$ and $\overset{\sim}{\boldsymbol\xi_0}$.
Therefore, the Dirichlet integral in (\ref{eq:Dir-Almgren})
is independent of the choice of orthonormal coordinates in $\mathbb{R}^n$.

Below, we recall from
Federer \cite{FH1} and Gr\"{u}ter \cite{Gr81}
some properties of functions in Sobolev spaces.
A map $X: \Omega\subset\mathbb{R}^{2}\rightarrow\mathbb{R}^N$
is called approximately differentiable at $w_0 \in\Omega$
with the approximate differential $\nabla X (w_0)$,
if there exists $X_0 \in \mathbb{R}^N$ such that for every $\epsilon>0$
\begin{equation*}
\Theta^2
\left(
\mathcal{L}^2 \lfloor
\Omega \smallsetminus
\{
w: |X(w)-X(w_0)-\nabla X(w_0) (w-w_0)| \le \epsilon \cdot |w-w_0|
\}
, w_0
\right)
=0 
\end{equation*}
where $\Theta^2$ denotes the two-dimensional density and
$\mathcal{L}^2 \lfloor D$
indicates the Lebesgue measure restricted to a set $D$
(see Federer {\cite[2.10.19]{FH1}}
or Gr\"{u}ter {\cite[Definition 2.2]{Gr81}}).  
Below is another characterization on the approximate differentiability. 
We say that 
$X: \Omega\subset\mathbb{R}^{2}\rightarrow\mathbb{R}^N$
is approximately differentiable at $w_0 \in\Omega$
with the approximate differential $\nabla X (w_0)$ 
if and only if there exists a neighborhood $U$ of $w_0$ and a map
$Y: U\rightarrow \mathbb{R}^N$ such that $Y$ is differentiable at $w_0$ and
\[
\Theta^2
\left(
\mathcal{L}^2 \lfloor
\Omega \smallsetminus
\{
w: X(w) \neq Y(w)
\}
, w_0
\right)
=0
\text{.}
\] 
The approximate differential is $\nabla Y (w_0)$.
If $X\in W^{1, 2}(\mathbb{R}^{2},\mathbb{R}^N)$, then
$X$ is approximately differentiable almost everywhere and the weak derivative
coincides with the approximate differential almost everywhere
(see Federer {\cite[Theorem 4.5.9 (26), (30) (VI))]{FH1}}.

\begin{definition}
[the set of ``good" points of a function in the Sobolev space $ W^{1, 2}(\Omega)$,
cf. Gr\"{u}ter \cite{Gr81}]
\label{def:good_point_sets_0}
Suppose $\Omega$ is a domain in $\mathbb{R}^2$ and 
$X\in W^{1, 2}(\Omega,\mathbb{R}^{N})$. 
Let 
$e(X)(w):=\left| \nabla X \right|^{2}(w)$
be the energy density of $X$. 
Define the set of ``good" points of
$X\in W^{1, 2}(\Omega,\mathbb{R}^{N})$
by
\begin{eqnarray*}
&A:=&
\{
w\in \Omega: X \text{ is approximately differentiable at $w$,}
\\
&&
\text{
$w$ is a Lebesgue point of $e(X)$,
$e(X)(w)
\neq 0$
}
\}
\text{.}
\end{eqnarray*}
\end{definition}

Below, we collect some lemmas from \cite{Gr81}.
We denote by
$\mathbb{U}_{r}(w_{0})\subset\mathbb{R}^{2}$ the open ball
$\{x\in\mathbb{R}^{2}: |x-w_{0}|<r \}$
and
$\mathbb{B}_{r}(w_{0})\subset\mathbb{R}^{2}$ the closed ball
$\{x\in\mathbb{R}^{2}: |x-w_{0}|\le r\}$.

\begin{lemma}[Gr\"{u}ter {\cite[2.5]{Gr81}}]
\label{lem:Gr2.5}
Let $X\in W^{1, 2}(\Omega,\mathbb{R}^{N})$ satisfy the
conformal conditions,
\begin{equation*}
\left| X_{u}\right| ^{2}=\left| X_{v}\right| ^{2},\text{ }X_{u}\cdot X_{v}=0
\text{, a.e. in $\Omega$.}
\end{equation*}
Suppose $\Omega$ is open and $w^{\ast }\in A\cap \Omega$.
Then,
\begin{equation*}
\underset{\sigma \rightarrow 0}{\lim \sup }\text{ }\sigma ^{-2}\text{ }
\underset{\Omega \cap \left\{ w:\left| X(w)-X(w^{\ast})\right|
<\sigma \right\} }{\int }\left| \nabla X \right|^{2}
\geq 2\pi
\text{.}
\end{equation*}
\end{lemma}

\begin{lemma}[Courant-Lebesgue Lemma, Gr\"{u}ter {\cite[2.6]{Gr81}}]
\label{lem:CL}
There is a constant $C=C(N)>0$ with the following property.
For any open set
$\Omega \subset \mathbb{R}^{2}$,
any $X\in W^{1, 2}(\Omega ,\mathbb{R}^{N})$,
any $w_{0}\in \Omega $,
and any $0<R<\text{dist} (w_{0},\partial \Omega )$,
there exists $r\in \lbrack \frac{1}{2}R,R]$ such that
\begin{equation*}
\underset{\partial \mathbb{U}_{r}(w_{0})}{\text{osc}} X
\leq C(N) \cdot
\left(\underset{\mathbb{U}_{R}(w_{0})}
{\int }\left| \nabla X \right|^{2}\right)^{1/2}
\text{.}
\end{equation*}
\end{lemma}

\section {The interior continuity}

In Definition \ref{def:WH},
we define the class of weakly harmonic multiple-valued functions
in the Sobolev space
$\mathcal{Y}_{2}(\Omega ,\mathbf{Q}_{Q}(\mathbb{R}^n))$.
The perturbations
in Definition \ref{def:WH}
are induced from the \texttt{range-variations},
which are also called outer variations in \cite{DS11}.

\begin{definition}[weakly harmonic multiple-valued functions]
\label{def:WH}
Let $\Omega\subset\mathbb{R}^{m}$ be an open set
and $\epsilon >0$ be sufficiently small.
For any given
$\psi\in C^{\infty}_{c}(\Omega\times\mathbb{R}^n,\mathbb{R}^n)$
such that the support
$\text{spt}(\psi)\subset\Omega^{\prime}\times\mathbb{R}^n$
for some
$\Omega^{\prime}\subset\subset\Omega$,
define the induced perturbation of $f$ by
\begin{equation}
f^{t}(x):= \sum\limits_{i=1}^{Q} [\![ f_{i}(x) + t \cdot \psi(x,f_{i}(x)) ]\!] 
\label{eq:Range_Variations}
\end{equation}
where $t\in (-\epsilon,\epsilon)$.
Then we say that
$f\in \mathcal{Y}_{2}(\Omega,\mathbf{Q}_{Q}(\mathbb{R}^n))$
is \texttt{weakly harmonic} if and only if
\begin{equation}
\lim_{t\rightarrow 0}\frac{Dir(f^{t};\Omega)-Dir(f;\Omega)}{t}=0
\text{.}
\label{eq:d_t(Dir(F_t-f))=0}
\end{equation}
\end{definition}

In Definition \ref{def:WNH},
we define the class of weakly Noether harmonic multiple-valued functions
in the Sobolev space
$\mathcal{Y}_{2}(\Omega ,\mathbf{Q}_{Q}(\mathbb{R}^n))$
of multiple-valued functions.
The perturbations
in Definition \ref{def:WNH}
are induced from the \texttt{domain-variations}, which are also called
inner variations in \cite{DS11}.

\begin{definition}[weakly Noether harmonic multiple-valued functions]
\label{def:WNH}
Let $\Omega\subset\mathbb{R}^{m}$ be an open set
and $\epsilon >0$ be sufficiently small.
Assume that, for any given $\phi\in C^{\infty}_{c}(\Omega,\mathbb{R}^m)$ and any fixed
$t\in(-\epsilon,\epsilon)$,
$X^{t}:\Omega\rightarrow\Omega$ defined by
$X^{t}(x):=x+t\cdot \phi(x)$
is a diffeomorphism of $\Omega$, leaving the boundary $\partial\Omega$ fixed.
We say that
$f\in \mathcal{Y}_{2}(\Omega,\mathbf{Q}_{Q}(\mathbb{R}^n))$
is \texttt{weakly Noether harmonic}
if and only if
\begin{equation*}
\lim_{t\rightarrow 0}\frac{Dir(f\circ X^{t};\Omega)-Dir(f;\Omega)}{t}=0
\text{.}
\end{equation*}
\end{definition}

\begin{definition}[weakly stationary-harmonic multiple-valued functions]
\label{def:WSH}
We say that
$f\in \mathcal{Y}_{2}(\Omega,\mathbf{Q}_{Q}(\mathbb{R}^n))$
is \texttt{weakly stationary-harmonic} if and only if
$f$ is \texttt{weakly harmonic} and \texttt{weakly Noether harmonic}.
\end{definition}

The main goal of this article is to prove the interior continuity
of any two-dimensional weakly stationary-harmonic multiple-valued function
in the Sobolev space
$\mathcal{Y}_{2}(\Omega ,\mathbf{Q}_{Q}(\mathbb{R}^n))$
of multiple-valued functions.

\begin{theorem}
\label{thm:main}
Suppose
$f\in \mathcal{Y}_{2}\left(\Omega, \mathbf{Q}_{Q}(\mathbb{R}^n)\right)$
is a strictly defined multiple-valued function,
where $\Omega$ is a simply connected open subset in $\mathbb{R}^{2}$.
Then $f$ is continuous in the interior of $\Omega $
as $f$ is weakly stationary-harmonic with finite Dirichlet integral over $\Omega$.
\end{theorem}

\subsection{The domain-variations:}

For convenience, we identify $\mathbb{C}$ with $\mathbb{R}^2$ below.
Note that Proposition \ref{prop:wc} is a well-known result
in the theory of $2$-dimensional harmonic mappings into Riemannian manifolds.

\begin{proposition}
\label{prop:wc}
Denote by
$\mathbb{U}_{R_0}(0)\subset\mathbb{C}$ the open ball of radius $R_0>0$
with center at the origin of complex plane $\mathbb{C}$.
Suppose
$f \in \mathcal{Y}_{2}(\mathbb{U}_{R_0}(0), \mathbf{Q}_{Q}(\mathbb{R}^n))$
is weakly stationary-harmonic.
Then,
\begin{enumerate}

\item[(1)]  The \texttt{Hopf differential} of
$\boldsymbol\xi_{0} \circ f
\in W^{1, 2}(\mathbb{U}_{R_0}(0), \mathbb{R}^{nQ})$,
\begin{equation*}
\begin{array}[t]{l}
\Phi (z):=
\left[ \left( \left| \frac{\partial \left(\boldsymbol\xi_{0} \circ f\right) }
{\partial u}\right| ^{2}-\left| \frac{\partial \left(\boldsymbol\xi_{0} \circ f\right)}
{\partial v}\right| ^{2}\right) -2i\text{ }\left\langle
\frac{\partial \left(\boldsymbol\xi_{0} \circ f\right) }{\partial u},
\frac{\partial \left(\boldsymbol\xi_{0} \circ f\right) }{\partial v}\right\rangle \right]
\text{ } d z^{2}
=: \varphi(z) \text{ } d z^{2} 
\end{array}
\end{equation*}
is holomorphic in the interior of $\mathbb{U}_{R_0}(0)$.
Here, $z=u+iv$ is a complex variable
and
$\left\langle \cdot , \cdot\right\rangle$
denotes the inner product of vectors in Euclidean spaces.

\item[(2)]  There exists a harmonic function
$h:\mathbb{U}_{R_0}(0)\rightarrow\mathbb{R}^2 \cong \mathbb{C}$
fulfilling
\begin{equation}
\left\{
\begin{array}{l}
\frac{\partial}{\partial z} h=\frac{-1}{4}\varphi
\text{,}
\\
\frac{\partial}{\partial \bar z}h=1
\text{,}
\end{array}
\right.
\label{eq:Harm_Fcn-0}
\end{equation}
such that
\begin{equation*}
(\boldsymbol\xi_{0} \circ f, h)
\in W^{1, 2}(\mathbb{U}_{R_0}(0), \mathbb{R}^{nQ}\times \mathbb{R}^2)
\end{equation*}
is weakly conformal in $\mathbb{U}_{R_0}(0)$, i.e.,
\begin{equation*}
\begin{array}[t]{l}
\left| \frac{\partial (\boldsymbol\xi_{0} \circ f, h)}{\partial u}\right|
=\left| \frac{\partial (\boldsymbol\xi_{0} \circ f, h)}{\partial v}\right| \text{ and }
\left\langle
\frac{\partial (\boldsymbol\xi_{0} \circ f, h)}{\partial u},
\frac{\partial (\boldsymbol\xi_{0} \circ f, h)}{\partial v}
\right\rangle
=0
\text{, a.e. in } \mathbb{U}_{R_0}(0)
\text{.}
\end{array}
\end{equation*}

\item[(3)]
Suppose $\boldsymbol\xi_{0}$, $\overset{\sim}{\boldsymbol\xi}_{0}$
are two distinct Lipschitz correspondences
and
$\varphi$, $\overset{\sim}{\varphi}$
are the induced holomorphic functions defined in (1) respectively.
Then,
$\lvert\varphi-\overset{\sim}{\varphi}\rvert
= \lvert\varphi_2-\overset{\sim}{\varphi}_2\rvert
=C(\boldsymbol\xi_{0}, \overset{\sim}{\boldsymbol\xi}_{0})$
is a constant depending only on the choice of
$\boldsymbol\xi_{0}$ and $\overset{\sim}{\boldsymbol\xi}_{0}$.
Moreover,
\begin{equation}
C(\boldsymbol\xi_{0}, \overset{\sim}{\boldsymbol\xi}_{0})
\le
\frac{4}{\pi R_0^2}\cdot Dir(f; \mathbb{U}_{R_0}(0))
\text{.}
\label{eq:Upper_Bdd(varphi-bar_varphi)}
\end{equation}
\end{enumerate}
\end{proposition}

\proof

$(1)$
We sketch the proof from \cite{Schoen84}.
For any smooth function $\eta:\mathbb{U}_{R_0}(0)\rightarrow\mathbb{R}$
with compact support,
let
$X^{t}(u,v)=(u+t\cdot \eta (u,v),v)$.
By the chain rule for weak derivatives, we have
\begin{equation*}
\left\{
\begin{array}{l}
\frac{\partial \left(\boldsymbol\xi_{0}\circ f \circ X^{t}(u,v)\right)}{\partial u}
= \frac{\partial \left(\boldsymbol\xi_{0}\circ f \right)}{\partial u}
\left( X^{t}(u,v) \right) \cdot\left(1+ t\cdot \frac{\partial \eta}{\partial u} \right)
\text{,}
\\ \\
\frac{\partial \left(\boldsymbol\xi_{0}\circ f \circ X^{t}(u,v)\right)}{\partial v}
= \frac{\partial \left(\boldsymbol\xi_{0}\circ f \right)}{\partial u}
\left(X^{t}(u,v)\right) \cdot\left( t\cdot \frac{\partial \eta}{\partial v} \right)
+
\frac{\partial \left(\boldsymbol\xi_{0}\circ f \right)}{\partial v}
\left( X^{t}(u,v) \right)
\text{.}
\end{array}
\right.
\end{equation*}
From the change of variables, $(\sigma,\tau)=X^{t}(u,v)$,
and the definition of Dirichlet integral of
multiple-valued functions in (\ref{eq:Dir-Almgren}), we have
\begin{equation*}
\begin{array}[t]{l}
Dir(f^{t};\Omega)=
\underset{\Omega}{\int} \text{ } \left[
\left|
\frac{\partial \left(\boldsymbol\xi_{0}\circ f \right)}{\partial \sigma}
\right|^{2}
\cdot
\left(1+ t \frac{\partial \eta}{\partial u} \right)^{2}
+
\left|
\frac{\partial \left(\boldsymbol\xi_{0}\circ f \right)}{\partial \sigma} \cdot
\left( t \frac{\partial \eta}{\partial v} \right)
+ \frac{\partial \left(\boldsymbol\xi_{0}\circ f \right)}{\partial \tau}
\right|^{2}
\right]
\text{ }
\frac{du dv}{1+t\cdot \eta_{u}}
\text{.}
\end{array}
\end{equation*}
Then, the vanishing of the first variations of the Dirichlet integral of $f$
with respect to the domain-variations $X^{t}$ gives
\begin{equation*}
\begin{array}[t]{l}
\int\limits_{\Omega}\text{ }
\left[ \left( \left| \frac{\partial \left(\boldsymbol\xi_{0} \circ f\right) }
{\partial u}\right| ^{2}-\left| \frac{\partial \left(\boldsymbol\xi_{0} \circ f\right)}
{\partial v}\right| ^{2}\right) \frac{\partial \eta }{\partial u}
+2
\left\langle
\frac{\partial \left(\boldsymbol\xi_{0} \circ f\right) }{\partial u},
\frac{\partial \left(\boldsymbol\xi_{0} \circ f\right) }{\partial v}
\right\rangle
\frac{\partial \eta }{\partial v}\right] \text{ }du dv=0
\text{.}
\end{array}
\end{equation*}
A similar argument, using the diffeomorphism
$X^{t}(u,v)=(u,v+t\cdot \zeta (u,v))$, gives
\begin{equation*}
\begin{array}[t]{l}
\int\limits_{\Omega}\text{ }
\left[
\left(
\left| \frac{\partial \left(\boldsymbol\xi_{0} \circ f\right) }
{\partial u}\right| ^{2}
-\left| \frac{\partial \left(\boldsymbol\xi_{0} \circ f\right)}
{\partial v}\right| ^{2}
\right)
\frac{\partial \zeta }{\partial v}
-2
\left\langle
\frac{\partial \left(\boldsymbol\xi_{0} \circ f\right) }{\partial u},
\frac{\partial \left(\boldsymbol\xi_{0} \circ f\right) }{\partial v}
\right\rangle
\frac{\partial \zeta }{\partial u}\right] \text{ } du dv=0
\text{.}
\end{array}
\end{equation*}
These two equations provide the weak form of the Cauchy-Riemann
equations for the $L^{1}$-function

\begin{equation}
\begin{array}[t]{l}
\varphi (z)=
\left(
\left|
\frac{\partial \left(\boldsymbol\xi_{0} \circ f\right) }{\partial u}
\right|^{2}
-
\left|
\frac{\partial \left(\boldsymbol\xi_{0} \circ f\right)}{\partial v}
\right|^{2}
\right)
-2i
\left\langle
\frac{\partial \left(\boldsymbol\xi_{0} \circ f\right) }{\partial u},
\frac{\partial \left(\boldsymbol\xi_{0} \circ f \right) }{\partial v}
\right\rangle
\text{.}
\end{array}
\label{eq:varphi=varphi_1+i.varphi_2}
\end{equation}
By Weyl's lemma, $\varphi $ is a holomorphic function of $z$.

$(2)$
If $\Phi$ is identically zero, then $\boldsymbol\xi_{0} \circ f$ is weakly
conformal.
From (\ref{eq:Harm_Fcn-0}), the assertion is then proved by choosing 
$h=\overline{z}$
(i.e., $h=u-i v$).
Therefore, we assume below that $\Phi$ is not identically zero.
Below we would like to follow the trick in \cite{Gr86}
to construct a $\mathbb{R}^2$-valued harmonic function
$h$ showing that the Hopf differential of
\begin{equation*}
(\boldsymbol\xi_{0} \circ f, h):\mathbb{U}_{R_0}(0)
\rightarrow \mathbb{R}^{nQ}\times \mathbb{R}^2
\end{equation*}
is identically zero.

For convenience, we abuse the notation by writing
$h=(h_{1}, h_{2})\in\mathbb{R}^2$ as
$h=h_{1}+i \cdot h_{2} \in\mathbb{C}$ below.
Let $\Phi_{h}(z)$ denote the Hopf differential of
$(\boldsymbol\xi_{0} \circ f, h)$.
By a simple computation, one can verify
\begin{equation}
\begin{array}[t]{l}
\Phi_{h} (z)
=\Phi (z)
+\left[
\left( \left| \frac{\partial h}{\partial u}\right|^{2}
-\left| \frac{\partial h}{\partial v}\right|^{2}\right)
-2i
\left\langle \frac{\partial h}{\partial u},
\frac{\partial h}{\partial v}
\right\rangle
\right] \text{ }dz^{2}
\\ \\
=
\left[
\varphi (z)+4 \text{ } \frac{\partial h}{\partial z} \text{ }
\frac{\partial \overline{h}}{\partial z}(z)\right]
\text{ }dz^{2}
\text{.}
\label{eq:wc-1}
\end{array}
\end{equation}
Since we know that $\varphi$ is holomorphic, there exists a holomorphic function
$\psi$ satisfying $\psi^{\prime}=\frac{-1}{4}\varphi$.
Let
\begin{equation}
h(z)=\psi (z)+\bar z
\text{.}
\label{eq:h=psi+bar(z)}
\end{equation}
Then, it is easy to verify that both
$h_{1}$ and $h_{2}$ are real-valued harmonic functions.
Moreover, a simple calculation shows that $h$ also fulfills
\begin{equation}
\begin{array}[t]{l}
\frac{\partial h}{\partial z} \text{ } \frac{\partial \overline{h}}{\partial z}
= -\frac{1}{4} \varphi
\text{.}
\end{array}
\label{eq:wc-2}
\end{equation}
From (\ref{eq:wc-1}) and (\ref{eq:wc-2}), $\Phi_{h} (z)\equiv 0$.
Therefore, one concludes that
$(\boldsymbol\xi_{0} \circ f, h): \mathbb{U}_{R_0}(0)
\rightarrow \mathbb{R}^{nQ}\times \mathbb{R}^2$
is weakly conformal.

$(3)$
Let
$\varphi=\varphi_1+i\cdot\varphi_2$ and
$\overset{\sim}{\varphi}
=\overset{\sim}{\varphi}_1+i\cdot\overset{\sim}{\varphi}_2$.
Notice that, from (\ref{eq:eq_xi_0-f}), we have
\begin{equation*}
\left\{
\begin{array}{l}
\lvert\partial_{u}(\boldsymbol\xi_{0}\circ f)\rvert
=
\lvert\partial_{u}(\overset{\sim}{\boldsymbol\xi}_{0}\circ f)\rvert
\text{,}
\\
\lvert\partial_{v}(\boldsymbol\xi_{0}\circ f)\rvert
=
\lvert\partial_{v}(\overset{\sim}{\boldsymbol\xi}_{0}\circ f)\rvert
\text{.}
\end{array}
\right.
\end{equation*}
From (\ref{eq:varphi=varphi_1+i.varphi_2}), 
$\overset{\sim}{\varphi}_1=\varphi_1$.
Since both $\varphi$ and $\overset{\sim}{\varphi}$ are holomorphic,
they satisfy the Cauchy-Riemann equations.
Thus,
\begin{equation*}
\left\{
\begin{array}{l}
\partial_{u}(\overset{\sim}{\varphi}_2-\varphi_2)
=-\partial_{v}(\overset{\sim}{\varphi}_1-\varphi_1)=0 \text{,}
\\
\partial_{v}(\overset{\sim}{\varphi}_2-\varphi_2)
=\partial_{u}(\overset{\sim}{\varphi}_1-\varphi_1)=0 \text{.}
\end{array}
\right.
\end{equation*}
Therefore, one concludes that
$\lvert\overset{\sim}{\varphi}-\varphi \rvert
= \lvert\overset{\sim}{\varphi}_2-\varphi_2 \rvert
= C(\boldsymbol\xi_{0}, \overset{\sim}{\boldsymbol\xi}_{0})$
is a constant depending on the choice of
$\boldsymbol\xi_{0}$ and $\overset{\sim}{\boldsymbol\xi}_{0}$
(or the choice of coordinates of $\mathbb{R}^{n}$).
Note that, from (\ref{eq:varphi=varphi_1+i.varphi_2}), we have
\begin{equation}
|\varphi|\le 2\cdot |\nabla(\boldsymbol\xi_{0}\circ f)|^2
\text{; }
|\overset{\sim}{\varphi}|\le 2\cdot |\nabla(\overset{\sim}{\boldsymbol\xi}_{0}\circ f)|^2
\text{.}
\label{eq:|varphi|<|nabla_f|^2}
\end{equation}
Thus, 
\begin{equation}
C(\boldsymbol\xi_{0}, \overset{\sim}{\boldsymbol\xi}_{0})
=|\overset{\sim}{\varphi}-\varphi|
\le
|\overset{\sim}{\varphi}|+|\varphi|
\le
2\cdot |\nabla(\boldsymbol\xi_{0}\circ f)|^2
+ 2\cdot |\nabla(\overset{\sim}{\boldsymbol\xi}_{0}\circ f)|^2
\text{.}
\label{eq:C(h,sim_h)}
\end{equation}
From integration over the set $\mathbb{U}_{R_0}(0)$,
(\ref{eq:C(h,sim_h)}) gives
\begin{equation*}
C(\boldsymbol\xi_{0}, \overset{\sim}{\boldsymbol\xi}_{0})
\le
\frac{4}{\pi R_0^2}\cdot
Dir(f; \mathbb{U}_{R_0}(0))
\text{.}
\end{equation*}

\endproof

\subsection{The range-variations:}

We would like to follow the proof of Theorem 3.10 in
Gr\"{u}ter's paper \cite{Gr81} to derive the key estimates.
In order to apply Gr\"{u}ter's argument,
we build up the so-called \emph{nested} \emph{admissible} closed balls of a member
$q\in\mathbf{Q}_{Q}(\mathbb{R}^n)$ 
(see Definition \ref{def:Nested_Adm_Closed_Ball})
for the construction of admissible range-variations.
For this purpose, we first need Proposition \ref{prop:Adm_Ball_Proj} and
Proposition \ref{prop:Nested_Adm_Closed_Ball}.

% below was modified on 09.05.2013
\begin{proposition}
\label{prop:Adm_Ball_Proj}
Suppose $n\ge 2$, $Q\ge 2$. 
Let $e_\alpha\in\{e_1,...,e_n, e_{n+1},...,e_{P(n,Q)}\}$ 
denote the unit tangent vector of the $\alpha$-th straight line in the bi-Lipschitz embedding 
$\boldsymbol\xi$, 
where $\{e_1,...,e_n\}$ forms the coordinate basis of $\mathbb{R}^n$. 
Then, for any
$q=\sum\limits_{j=1}^{Q} [\![q_j]\!]\in\mathbf{Q}_{Q}(\mathbb{R}^n)$, 
there exists a positive number 
\[
0<\theta_0=\theta_0(n,Q)
\le \frac{\pi}{4}
\] 
such that the inequality
\begin{equation}
\left|
\left<
e_{\alpha}, \frac{q_i-q_j}{|q_i-q_j|}
\right>
\right|
\ge \sin\theta_0
\label{eq:measuredangle_Lower_Bdd}
\end{equation}
holds for any $\alpha\in\{1,..., P(n,Q)\}$
and any $q_{i}\ne q_{j}$. 
\end{proposition}
% above was modified on 09.05.2013

\proof 
%below was added on 09.05.2013
Step $1^{\circ}$ 
%above was added on 09.05.2013 
Denote by $\pi_{\alpha}$ the hyperplane with normal vector 
$\pm e_{\alpha}$ passing through the origin of $\mathbb{R}^n$, and by 
$P_{\pi_\alpha}(v)$ the orthogonal projection of the vector $v\in \mathbb{R}^n$ 
into the hyperplane 
$\pi_{\alpha}$. 
Denote by 
$\measuredangle(v,\pi_{\alpha})$ the positive angle between the vector 
$v\in \mathbb{R}^n$ and the vector 
$P_{\pi_\alpha}(v)\in\pi_{\alpha}\subset\mathbb{R}^{n}$. 
We let
$\measuredangle(v, \pi_{\alpha})\in[0,\frac{\pi}{2}]$, 
and let $\measuredangle(v,\pi_{\alpha})=\frac{\pi}{2}$
as $P_{\pi_\alpha}(v)=0$. 
Denote by $\mathbb{S}^{n-1}$ the unit sphere in $\mathbb{R}^{n}$ 
with center at the origin of $\mathbb{R}^n$. 
As $v\in\mathbb{S}^{n-1}$, 
the distance on the unit sphere
between the point $v\in\mathbb{S}^{n-1}$ 
and the ``great circle" $\pi_{\alpha}\cap\mathbb{S}^{n-1}$ 
is meant to be the positive angle between the unit vector $v\in\mathbb{R}^{n}$ 
and the vector $P_{\pi_\alpha}(v)$. 
Let 
\[
\{v_{1},...,v_{\ell}\}:=\bigcup\limits_{i,j \in\{1,...,Q\}}
\left\{\frac{q_{i}-q_{j}}{|q_{i}-q_{j}|}:q_{i}\neq q_{j} \right\}
\]
for some 
\begin{equation}
\ell\le Q(Q-1) 
\text{,} 
\label{eq:ell<=Q(Q-1)}
\end{equation}
and 
$v_{i}\neq v_{j}$ for any distinct $i,j\in\{1,...,\ell\}$.

The proof is equivalent to showing that
there exist a positive number $\theta_{0}=\theta_{0}(n,Q)>0$
and a coordinate bases $\{e_1,...,e_n\}$ of $\mathbb{R}^{n}$
such that the distance on $\mathbb{S}^{n-1}$
between any member in the given set 
$\bigcup\limits_{i,j \in\{1,...,Q\}}$ 
$\left\{\frac{q_{i}-q_{j}}{|q_{i}-q_{j}|}:q_{i}\neq q_{j}\right\}$ 
$\subset\mathbb{S}^{n-1}$
and any member in the set 
$\bigcup\limits_{\alpha=1}^{n}$ $\pi_{\alpha}\cap\mathbb{S}^{n-1}$
is at least $\theta_{0}(n,Q)$.
In other words,
we need to show that we could choose a new coordinate bases 
$\{\overset{\circ}{e}_1,..., \overset{\circ}{e}_n\}$ 
of $\mathbb{R}^{n}$ such that
$\measuredangle\left(v_{i,j}, 
\overset{\circ}{\pi}_{\alpha}\right)\ge\theta_{0}(n,Q)>0$, 
$\forall$ $v_{i,j}\ne 0$ and $\forall$ $\alpha\in\{1,...,n\}$. 
Here, $\overset{\circ}{\pi}_{\alpha}$ denotes the hyperplane with normal vector 
$\pm \overset{\circ}{e}_{\alpha}$. 
Let
\begin{equation}
\delta_{\ell}:=\sin^{-1}\left(\frac{1}{\sqrt{n}}\right)
\cdot\left(\frac{1}{2}\right)^{(n-1)(\ell-1)}
\text{.}
\label{eq:def-delta_ell}
\end{equation}
In other words, we let 
$\delta_{1}=\sin^{-1}\left(\frac{1}{\sqrt{n}}\right)$ 
and 
$\delta_{\ell+1}=\delta_{\ell}\cdot\left(\frac{1}{2}\right)^{n-1}$,  
$\forall$ $\ell\in\mathbb{Z}_{+}$.

The proof is an induction argument. 
For any given point $v_{1}\in\mathbb{S}^{n-1}$,
we choose a new coordinate bases $\{e_{1}^{1},...,e_{n}^{1}\}$ of
$\mathbb{R}^{n}$ by rotating the original one so that
$v_{1}=\sum\limits_{\alpha=1}^{n} \frac{1}{\sqrt{n}}\cdot e_{\alpha}^{1}$. 
In other words,
$v_{1}=\left(\frac{1}{\sqrt{n}},...,\frac{1}{\sqrt{n}}\right)$
under the new coordinates, and
$\measuredangle\left(v_{1},\pi_{\alpha}\right)=\sin^{-1}\left(\frac{1}{\sqrt{n}}\right)$
for all $\alpha\in\{1,...,n\}$.

Let $\ell\ge2$ and
assume that, 
for any $(\ell-1)$ distinct points
$v_{1},...,v_{\ell-1}\in\mathbb{S}^{n-1}$, 
there exists a coordinate bases 
$\left\{e_{1}^{\ell-1},...,e_{n}^{\ell-1}\right\}$ of $\mathbb{R}^n$ 
such that
\begin{equation}
\measuredangle \left(v_{b},\pi_{\alpha}^{\ell-1}\right)
\ge \delta_{\ell-1}
\text{, }
\forall \text{ } b\in\{1,...,\ell-1\} \text{ and }\forall \text{ }\alpha\in\{1,...,n\} 
%\text{,}
\label{eq:The-(ell-1)-th-condition}
\end{equation}
where $\pi_{\alpha}^{\ell-1}$ is the hyperplane in $\mathbb{R}^n$ 
with normal vector
$\pm e_{\alpha}^{\ell-1}$. 
Then,
for any extra given point $v_{\ell}\in\mathbb{S}^{n-1}$, 
we will show that there exists a new coordinate bases 
$\left\{e_{1}^{\ell},...,e_{n}^{\ell}\right\}$ of $\mathbb{R}^n$ 
such that
\begin{equation}
\measuredangle \left(v_{b},\pi_{\alpha}^{\ell}\right)
\ge \delta_{\ell}
\text{, }
\forall \text{ } b\in\{1,...,\ell\} \text{ and }\forall \text{ }\alpha\in\{1,...,n\}
\text{.}
\label{eq:The-(ell)-th-condition}
\end{equation}
We may assume that, for the extra vector $v_{\ell}\ne 0$, 
$\exists$ $\alpha\in\{1,...,n\}$ such that 
\begin{equation}
\measuredangle\left(v_{\ell},\pi_{\alpha}^{\ell-1}\right)<\delta_{\ell-1} 
\text{.}
\label{eq:ell-th-angle-small}
\end{equation}
Otherwise, 
$\measuredangle\left(v_{\ell},\pi_{\alpha}^{\ell-1}\right)
\ge\delta_{\ell-1}>\delta_{\ell}, 
\text{ }\forall \alpha\in\{1,...,n\}$. 
Then from 
(\ref{eq:The-(ell-1)-th-condition}), 
we now obtain (\ref{eq:The-(ell)-th-condition}) by letting 
$e_{\alpha}^{\ell}:=e_{\alpha}^{\ell-1}$, $\forall$ $\alpha\in\{1,...,n\}$. 
On the other hand, (\ref{eq:ell-th-angle-small}) can't hold
for all hyperplane $\pi_{\alpha}^{\ell-1}$, $\alpha=1,...,n$.
Otherwise, from $\delta_{\ell}\le\sin^{-1}(\frac{1}{\sqrt{n}})$, 
it imples $|v_{\ell}|< 1$, 
which contradicts $v_{\ell}\in\mathbb{S}^{n-1}$.
Therefore, we now assume 
\begin{equation}
\measuredangle\left(v_{\ell},\pi_{\alpha_{1}}^{\ell-1}\right)
\le\cdots\le
\measuredangle\left(v_{\ell},\pi_{\alpha_{K}}^{\ell-1}\right)
<\delta_{\ell-1}\le
\measuredangle\left(v_{\ell},\pi_{\alpha_{K+1}}^{\ell-1}\right)
\le\cdots\le
\measuredangle\left(v_{\ell},\pi_{\alpha_{n}}^{\ell-1}\right)
\label{eq:moving_seq}
\end{equation}
for some integer $K$ fulfilling $1\le K\le n-1$.

The new coordinate bases $\{e_{1}^{\ell},...,e_{n}^{\ell}\}$ of $\mathbb{R}^{n}$ is obtained from rotation. 
Below, we construct a rotation of coordinate bases of $\mathbb{R}^{n}$ 
by considering it as how to move  
the extra point $v_{\ell}$ on the unit sphere. 
Denote by 
$v_{\ell}^{-}(i)$ the position of $v_{\ell}$ 
before moving $v_{\ell}$  
along the shortest geodesic connecting 
$v_{\ell}$ and $\pm e_{\alpha_{i}}^{\ell-1}$; 
and by $v_{\ell}^{+}(i)$ the new position of $v_{\ell}$ 
after this movement of $v_{\ell}$. 
We let $v_{\ell}^{+}(i)=v_{\ell}^{-}(i+1)$ for each $i\in\{1,...,K-1\}$.
Here, $\pm e_{\alpha_{1}}^{\ell-1}$ is chosen to be either 
$e_{\alpha_{1}}^{\ell-1}$ or $-e_{\alpha_{1}}^{\ell-1}$, 
depending on whether 
$\measuredangle\left(v_{\ell}, e_{\alpha_{1}}^{\ell-1}\right)$ 
or 
$\measuredangle\left(v_{\ell}, -e_{\alpha_{1}}^{\ell-1}\right)$ 
is the smallest one. 
Under the assumption of (\ref{eq:moving_seq}), 
we first move $v_{\ell}$ on the unit sphere along the shortest geodesic 
connecting $v_{\ell}$ and $\pm e_{\alpha_{1}}^{\ell-1}$ 
until
$\measuredangle\left(v_{\ell},\pi_{\alpha_{1}}^{\ell-1}\right)
=\frac{1}{2}\delta_{\ell-1}$. 
We proceed similarly for each $i\in\{2,...,K\}$ 
by moving $v_{\ell}^{-}(i)$ on the unit sphere along the shortest geodesic 
connecting $v_{\ell}^{-}(i)$ and $\pm e_{\alpha_{1}}^{\ell-1}$ 
until it arrives the position $v_{\ell}^{+}(i)$ which satisfies 
$\measuredangle\left(v_{\ell}^{+}(i),\pi_{\alpha_{i}}^{\ell-1}\right) 
=\left(\frac{1}{2}\right)^{i}\delta_{\ell-1}$. 
We let 
$v_{\ell}^{+}:=v_{\ell}^{+}(K)$ and $v_{\ell}^{-}:=v_{\ell}^{-}(1)$, 
and note 
\begin{equation}
\measuredangle\left(v_{\ell}^{-}(i),v_{\ell}^{+}(i)\right)
\le
\delta_{\ell-1}\cdot
\left(\frac{1}{2}\right)^{i}
\text{, } \forall \text{ } i\in\{1,...,K\} 
\text{,}
\label{eq:movement_on_sphere-1}
\end{equation}
during this procedure.

Now we describe the procedure of moving the point 
$v_{\ell}\in\mathbb{S}^{n-1}$ 
above by the formulation of rotations. 
Denote by $\mathcal{R}_{\ell}^{i}$ the rotation of coordinate bases of 
$\mathbb{R}^n$ corresponding to the movement, 
$v_{\ell}^{-}(i) \rightarrow v_{\ell}^{+}(i)$. 
Let 
\begin{equation}
\mathcal{R}_{\ell}:=\mathcal{R}_{\ell}^{K} 
\centerdot \centerdot \centerdot
\mathcal{R}_{\ell}^{1}
\label{eq:rotation_composition}
\end{equation}
be the rotation of coordinate bases of 
$\mathbb{R}^n$ corresponding to the movement, 
$v_{\ell}^{-}\rightarrow v_{\ell}^{+}$. 
We may associate the ``dual" rotation operator 
$\overset{\sim}{\mathcal{R}}$ 
acting on vectors in 
$\mathbb{R}^n$ by letting 
\begin{equation}
\measuredangle(
\overset{\sim}{\mathcal{R}}(v), \pi_{\alpha}^{\ell-1}) 
:= 
\measuredangle(v, \mathcal{R}(\pi_{\alpha}^{\ell-1}) ) 
\text{, }
\forall \text{ } \alpha\in\{1,...,n\}
\text{, } \forall \text{ } v\in \mathbb{R}^n 
\text{.}
\label{eq:rotation_duality}
\end{equation}
Now we let 
$\{ 
e_{1}^{\ell}, 
..., 
e_{n}^{\ell} 
\}$ 
be the coordinate bases of $\mathbb{R}^n$ 
derived from 
$e_{\alpha}^{\ell}:=\overset{\sim}{\mathcal{R}}_{\ell}(e_{\alpha}^{\ell-1})$, 
$\forall$ $\alpha\in\{1,...,n\}$.  
Denote by $\pi_{\alpha}^{\ell}:=\mathcal{R}_{\ell}(\pi_{\alpha}^{\ell-1})$ 
the hyperplane in $\mathbb{R}^n$ 
with the normal vector 
$\pm e_{\alpha}^{\ell}
:=\overset{\sim}{\mathcal{R}}_{\ell}(\pm e_{\alpha}^{\ell-1})$. 
Note that we have 
\begin{equation}
\measuredangle(v_{\ell}^{+}, \pi_{\alpha}^{\ell-1}) 
= 
\measuredangle(v_{\ell}^{-}, \mathcal{R}_{\ell}(\pi_{\alpha}^{\ell-1}) ) 
=
\measuredangle(v_{\ell}, \pi_{\alpha}^{\ell}) 
\text{, }
\forall \text{ } \alpha\in\{1,...,n\}
\text{.}
\label{eq:moving_equiv_to_rotation}
\end{equation} 
Moreover, according to (\ref{eq:movement_on_sphere-1}), 
the angle of rotation $\mathcal{R}_{\ell}^{i}$ 
(or $\overset{\sim}{\mathcal{R}_{\ell}^{i}}$) is no greater than 
$\delta_{\ell-1}\left(\frac{1}{2}\right)^{i}$. 
This implies that, 
$\forall$ $b \in\{1,...,\ell\}$ and $\forall$ $\alpha\in\{1,...,n\}$, 
\begin{equation}
\measuredangle(v_{b}, \pi_{\alpha}^{\ell-1}) 
-\delta_{\ell-1}\left(\frac{1}{2}\right)^{i}
\le 
\measuredangle(v_{b}, \mathcal{R}_{\ell}^{i}(\pi_{\alpha}^{\ell-1}) ) 
\le 
\measuredangle(v_{b}, \pi_{\alpha}^{\ell-1}) 
+\delta_{\ell-1}\left(\frac{1}{2}\right)^{i}
\text{.}
\label{eq:rotation_angle}
\end{equation} 
To see it, observe that 
\[
\measuredangle(v_{b}, \mathcal{R}_{\ell}^{i}(\pi_{\alpha}^{\ell-1}) ) 
= \measuredangle(\overset{\sim}{\mathcal{R}_{\ell}^{i}}(v_{b}), 
\pi_{\alpha}^{\ell-1}) 
\le 
\measuredangle(\overset{\sim}{\mathcal{R}_{\ell}^{i}}(v_{b}), 
v_{b}) 
+
\measuredangle(v_{b}, 
\pi_{\alpha}^{\ell-1}) 
\text{.}
\] 
From (\ref{eq:movement_on_sphere-1}), 
we obtain the second inequality in 
(\ref{eq:rotation_angle}). 
Similarly, by letting 
$\mathcal{S}_{\ell}^{i}:=(\mathcal{R}_{\ell}^{i})^{-1}$, 
the inverse of $\mathcal{R}_{\ell}^{i}$, 
and observe from 
\[
\measuredangle(v_{b}, \pi_{\alpha}^{\ell-1})
=\measuredangle(
\overset{\sim}{\mathcal{S}_{\ell}^{i}}(v_{b}), 
\mathcal{R}_{\ell}^{i}(\pi_{\alpha}^{\ell-1}) ) 
\le 
\measuredangle(\overset{\sim}{\mathcal{S}_{\ell}^{i}}(v_{b}), 
v_{b}) 
+
\measuredangle(v_{b}, 
\mathcal{R}_{\ell}^{i}(\pi_{\alpha}^{\ell-1}) ) 
\text{,}
\]  
we obtain the first inequality in 
(\ref{eq:rotation_angle}). 
Notice that we have applied the fact that 
$\measuredangle(\overset{\sim}{\mathcal{S}_{\ell}^{i}}(v_{b}), 
v_{b}) = \measuredangle(\overset{\sim}{\mathcal{R}_{\ell}^{i}}(v_{b}), 
v_{b}) $. 

From (\ref{eq:rotation_angle}), (\ref{eq:moving_equiv_to_rotation}), 
(\ref{eq:rotation_composition}) and $K\le n-1$,  
we conclude 
\[
\measuredangle(v_{b}, \pi_{\alpha}^{\ell}) 
\ge 
\delta_{\ell-1}\cdot
\left( \frac{1}{2} \right)^{n-1}
\text{, }\forall\text{ } \alpha\in\{1,...,n\} 
\text{, } \forall\text{ } b\in\{1,...,\ell\}
\text{.}
\] 
Now, from (\ref{eq:def-delta_ell}) and (\ref{eq:ell<=Q(Q-1)}), 
we finish the proof as $\alpha\in\{1,2,...,n\}$ by letting 
% below was modified on 09.05.2013 
$
\theta_0=
\left(\frac{1}{2}\right)^{(n-1)(Q(Q-1)-1)}
\sin^{-1}\left(\frac{1}{\sqrt{n}}\right) 
\text{.}
$

Step $2^{\circ}$ 
For the cases of $\alpha\in\{n+1,...,P(n,Q)\}$, we attach the extra straight lines $\{L_{n+1},...,L_{P(n,Q)}\}$ 
to the coordinate basis $\{e_{1},...,e_{n}\}$ rigidly. 
Then, we proceed the argument above further in finite steps (at most $P(n,Q)-n$ times) and obtain a smaller 
$\theta_{0}(n,Q)\in(0,\frac{\pi}{4}]$ than the one in Step $1^{\circ}$.  
% above was modified on 09.05.2013 
\endproof

\begin{definition}
[The admissible closed balls of $q$ in $\mathbf{Q}_{Q}(\mathbb{R}^{n})$ and
the union of admissible closed balls of $\text{spt}(q)$ in $\mathbb{R}^{n}$]
\label{def:Adm_Closed_Ball}
%\label{def:Nested_Adm_Closed_Ball}
Let
$q=\sum\limits_{i=1}^{I} \ell_i [\![q_i]\!]$
be any member in $\mathbf{Q}_{Q}(\mathbb{R}^{n})$,
where $q_1,..., q_I$ are distinct points in $\mathbb{R}^{n}$.
Denote by
\begin{equation*}
\mathbb{B}^{n}_{\tau}\left(y_0\right)
:=\{y\in\mathbb{R}^{n}: |y-y_0| \le \tau \}
\end{equation*}
the closed ball of radius $\tau$ with center $y_0$ in $\mathbb{R}^{n}$, 
by
\begin{equation*}
\mathcal{B}_{\tau} (q)
:=\bigcup\limits_{i=1}^{I}\mathbb{B}^{n}_{\tau}\left( q_i \right)
\end{equation*}
the union of closed balls of radius
$\tau$ in $\mathbb{R}^{n}$,
and by
\begin{equation*}
\mathbb{B}^{\mathbf{Q}}_{\tau}(q):=
\left\{p\in\mathbf{Q}_{Q}(\mathbb{R}^{n}): \mathcal{G}(p, q)
\le \tau \right\}
\subset \mathbf{Q}_{Q}(\mathbb{R}^n) 
\end{equation*}
the closed ball of radius $\tau$ with center $q$ in
$\mathbf{Q}_{Q}(\mathbb{R}^{n})$.
Then,
$\mathcal{B}_{\tau} (q)$
is said to be \emph{admissible}
if 
\begin{equation*}
\Pi_{\varkappa}\left(\mathbb{B}_{\tau}^{n}(q_{i})\right)
\cap
\Pi_{\varkappa}\left(\mathbb{B}_{\tau}^{n}(q_{j})\right)
=\emptyset
\text{, }
\forall \text{ } \varkappa \in\{1,...,n\}
\text{, }
\forall \text{ } i\neq j \in\{1,...,I\}
\text{.}
\end{equation*}
The closed ball
$\mathbb{B}^{\mathbf{Q}}_{\tau}(q)$
is said to be \emph{admissible} if
$\mathcal{B}_{\tau} (q)$ is \emph{admissible}.

\end{definition}

From Proposition \ref{prop:Adm_Ball_Proj},
if
$\tau$ satisfies
\begin{equation*}
0<\tau<\frac{\sin\theta_0 (n,Q)}{2}\cdot\inf\limits_{i\ne j} \{ |q_{i}-q_{j}| \}
\end{equation*}
where $\theta_{0} (n,Q)$ is given in
Proposition \ref{prop:Adm_Ball_Proj} as $n\ge 2$
and $\theta_{0} (1,Q):=\pi/2$,
then $\mathcal{B}_{\tau} (q)$
is a union of admissible closed ball of radius $\tau$ in $\mathbb{R}^{n}$.

We introduce the notion of \emph{admissible} closed balls in
Definition \ref{def:Adm_Closed_Ball} for
the construction of smooth vector fields in $\mathbb{R}^n$
in the range-variations of multiple-valued functions.
There is a formula of ``subtraction" between
$q\in\mathbf{Q}_{Q}(\mathbb{R}^n)$
and any member in the admissible ball
$\mathbb{B}^{\mathbf{Q}}_{\tau}(q) \subset \mathbf{Q}_{Q}(\mathbb{R}^n)$
of $q$ (under a suitable choice of coordinate bases of $\mathbb{R}^n$).
In other words, for any
$p=\sum\limits_{j=1}^{Q} [\![p_j]\!]
\in\mathbb{B}^{\mathbf{Q}}_{\tau}(q)$,
an admissible closed ball of $q$ ($p_1,..., p_Q$ are not necessarily distinct),
there is a natural way to obtain an member in
$\mathbf{Q}_{Q}(\mathbb{R}^n)$,
denoted as $p\ominus q $ or $q\ominus p $,
such that
$\mathcal{G}\left( p\ominus q, Q[\![0]\!] \right)
=\mathcal{G}\left( p, q \right)$.
To explain it, observe that 
$p\in\mathbb{B}^{\mathbf{Q}}_{\tau}(q)$
implies
$\text{spt}(p)\subset \mathcal{B}_{\tau}(q)$
and
$\text{card} \left( \text{spt}(p)
\cap \mathbb{B}_{\tau}\left( q_i \right)\right) =\ell_i$
for each $i\in\{1,...,I\}$.
Otherwise, there exists $i\in\{1,...,I\}$ such that
$\text{card} \left( \text{spt}(p)
\cap \mathbb{B}_{\tau}\left( q_i \right)\right) \lvertneqq \ell_{i}$.
But this means that there is a point
$p_{\iota}\in \text{spt}(p)$ such that $|p_{\iota}-q_{i}|\gneqq \tau$,
which contradicts the assumption of
$p\in\mathbb{B}^{\mathbf{Q}}_{\tau}(q)$.
Thus, we may write
$p= \sum\limits_{j=1}^{Q} [\![p_j^{\varkappa_j}]\!]$,
where
$\text{card}\left(\left\{j\in\{1,...,Q\}: \varkappa_j=i \right\}\right)=\ell_{i}$,
$p_{j}^{\varkappa_j}\in \mathbb{B}_{\tau}\left( q_{\varkappa_j} \right)$,
$\varkappa_{j}\in\{1,...,I\}$,
for each $j\in\{1,...,Q\}$.
Furthermore,
for any $p$ in an admissible closed ball of $q$,
$\mathbb{B}_{\tau}^{Q}(q)$,
we can define the
``subtraction" $q\ominus p \in \mathbf{Q}_{Q}(\mathbb{R}^n)$ by
\begin{equation}
q\ominus p:= \sum\limits_{j=1}^{Q}
[\![ q_{\varkappa_j}- p_{j}^{\varkappa_j} ]\!]
\text{.}
\label{eq:def(p-q)}
\end{equation}
It is obvious that
\begin{equation}
\mathcal{G}\left( q\ominus p, Q[\![0]\!] \right)
=
\mathcal{G}\left( p, q \right)
\text{, }\forall p\in\mathbb{B}_{\tau}^{Q}(q)
\text{.}
\label{eq:G(p,g)=G(p-q,0)}
\end{equation}

Let
$q=\sum\limits_{i=1}^{I} \ell_{i} [\![ q_{i} ]\!]
\in\mathbf{Q}_{Q}(\mathbb{R}^{n})$.
Recall from (\ref{eq:Xi_0}) that we can choose the map
${\xi}(\Pi_{\alpha},\cdot ):\mathbf{Q}_{Q}(\mathbb{R}^{n})
\rightarrow\mathbb{R}^{Q}$ for all $\alpha\in\{1,...,n\}$.
Then, for each fixed $\alpha\in\{1,...,n\}$,
we may write
\begin{equation*}
{\xi}(\Pi_{\alpha}, q)=
\left\{
\left(
\Pi_{\alpha}\left( q_{\varsigma_{\alpha}(1)} \right),
\cdots,
\Pi_{\alpha}\left( q_{\varsigma_{\alpha}(I)} \right)
\right)
:
\Pi_{\alpha}\left( q_{\varsigma_{\alpha}(1)} \right)
\le \cdots \le \Pi_{\alpha}\left( q_{\varsigma_{\alpha}(I)} \right)
\right\}
\subset \mathbb{R}^{Q}
\end{equation*}
and
\begin{equation*}
{\xi}(\Pi_{\alpha}, p)=
\left\{
\left(
\Pi_{\alpha}\left( p_{\omega_{\alpha}(1)}^{\varsigma_{\alpha}(1)} \right),
\cdots,
\Pi_{\alpha}\left( p_{\omega_{\alpha}(Q)}^{\varsigma_{\alpha}(I)} \right)
\right)
:
\Pi_{\alpha}\left( p_{\omega_{\alpha}(1)}^{\varsigma_{\alpha}(1)} \right)
\le \cdots \le
\Pi_{\alpha}\left( p_{\omega_{\alpha}(Q)}^{\varsigma_{\alpha}(I)} \right)
\right\}
\end{equation*}
for some permutation
$\varsigma_\alpha:\{1,...,I\}\rightarrow\{1,...,I\}$
and
$\omega_\alpha:\{1,...,Q\}\rightarrow\{1,...,Q\}$.
Notice that,
as a member $p\in\mathbf{Q}_{Q}(\mathbb{R}^n)$ is contained in
an admissible closed ball of $q$,
the two types of ``subtraction",
$q\ominus p$
and
${\xi}(\Pi_{\alpha}, q)-{\xi}(\Pi_{\alpha}, p) \in\mathbb{R}^Q$ 
are consistent, 
$\forall \text{ } \alpha\in\{1,...,n\}$.
In other words, if
$\mathbb{B}^{\mathbf{Q}}_{\tau}(q)
\subset\mathbf{Q}_{Q}(\mathbb{R}^n)$
is an admissible closed ball of $q$ and
$p\in\mathbb{B}^{\mathbf{Q}}_{\tau}(q)$,
then we may define
$p(s):=s \cdot q+(1-s)\cdot p$ and $s\in\mathbb{R}$
by
\begin{equation}
p(s)
= \sum\limits_{j=1}^{Q} \text{ }
[\![
s\cdot \left(q_{\varkappa_j}- p_{j}^{\varkappa_j}\right)
+(1-s)\cdot p_{j}^{\varkappa_j}
]\!] 
\label{eq:The_subtraction}
\end{equation}
and $p(s)$ satisfies the property
\begin{equation*}
{\xi}\left(\Pi_{\alpha}, p(s) \right) =
s \cdot {\xi}(\Pi_{\alpha}, p(1)) + (1-s)\cdot {\xi}(\Pi_{\alpha}, p(0))
\text{.}
%\label{eq:The_linear_relation-1}
\end{equation*}
Thus, 
\begin{equation}
\boldsymbol\xi_{0} \left(p(s) \right) =
s \cdot \boldsymbol\xi_{0}(q) + (1-s) \cdot \boldsymbol\xi_{0}(p)
\text{,}
\label{eq:The_linear_relation}
\end{equation}
for all $p\in\mathbb{B}^{\mathbf{Q}}_{\tau}(q)$
(an admissible closed ball of $q$).
From
(\ref{eq:def(p-q)}), (\ref{eq:G(p,g)=G(p-q,0)})
and taking
$\frac{d}{ds}$ in (\ref{eq:The_linear_relation}), 
we conclude that
\begin{equation}
|\boldsymbol\xi_{0}(q\ominus p)|^2
=\left| \boldsymbol\xi_{0}(q)-\boldsymbol\xi_{0}(p) \right|^{2}
=\left[\mathcal{G}\left( q\ominus p, Q[\![0]\!] \right)\right]^{2}
=\left[\mathcal{G}\left( p, q \right)\right]^{2}
\text{,}
\label{eq:|xi(p)-xi(q)|=G(p,q)}
\end{equation}
for any $p$ in an admissible closed ball $\mathbb{B}^{\mathbf{Q}}_{\tau}(q)$.

\begin{proposition}
\label{prop:Nested_Adm_Closed_Ball}
For any
$q=\sum\limits_{i=1}^{I} \ell_i [\![q_i]\!]\in\mathbf{Q}_{Q}(\mathbb{R}^n)$,
there correspond a non-negative integer $L\in\{0,1,...,Q-1\}$,
a sequence of members
$q^{(0)}=q, q^{(1)},...,q^{(L)}\in \mathbf{Q}_{Q}(\mathbb{R}^n)$,
where
$q^{(k)}=\sum\limits_{i=1}^{I_k} \ell^k_i [\![q_i^k]\!]$,
a constant $C_0=C_0 (n,Q)>1$,
and a sequence of numbers,
\begin{equation}
0=:\rho_0<\sigma_0<\rho_1<\sigma_1<\rho_2<\cdots<\sigma_{L-1}
<\rho_L< \sigma_L:=+\infty
\label{eq:(rho-sigma)-sequence}
\end{equation}
such that the following statements hold.

\begin{enumerate}
\item[(a)]
If
$\text{card}\left(\text{spt}(q)\right) =1$,
then
$L= 0$,
$\rho_0=0$ and $\sigma_0=+\infty$.

\item[(b)]
If
$\text{card}\left(\text{spt}(q)\right) \ge 2$,
then
$L\ge 1$,
\begin{equation}
\text{spt}\left(q^{(k-1)}\right)
\supsetneq \text{spt}\left(q^{(k)}\right)
\label{eq:spt{q^k}-decreasing}
\end{equation}
for each $k\in\{1,...,L\}$ and
the set $\text{spt}\left(q^{(L)}\right)$
consists of a single point.

\item[(c)]
For each $k\in\{0,1,...,L-1\}$ and $L\ge 1$,
\begin{equation}
10 Q\cdot \rho_{k} \le \sigma_{k} 
\label{eq:Lower_Bdd-sigma_(k)/rho_(k)}
\end{equation}
where
\begin{equation}
\sigma_k:=\frac{\sin\theta_{0}}{4}\cdot\inf\limits_{i\ne j}
\left\{|q_i^k-q_j^k|)\right\} 
\label{eq:sigma_(k)-def}
\end{equation}
$\theta_0 (n,Q)\in(0,\pi/4)$ is given in
Proposition \ref{prop:Adm_Ball_Proj} as $n\ge2$
and $\theta_0 (n,Q)=\pi/2$ as $n=1$.

\item[(d)]
For each $k\in \{1,...,L\}$ and $L\ge 1$,
\begin{equation}
\sigma_{k-1}<\rho_{k} \le C_0(n,Q) \cdot \sigma_{k-1}
\text{.}
\label{eq:Upper_Bdd-rho_(k+1)/sigma_(k)}
\end{equation}

\item[(e)]
For each $k\in\{1,...,L\}$ and $L\ge 1$,
if $\sigma_{k-1}<\rho_{k}$,
then
\begin{equation}
\mathbb{B}^{\mathbf{Q}}_{\sigma_{k-1}}\left(q^{(k-1)}\right)
\subset\mathbb{B}^{\mathbf{Q}}_{\rho_{k}}\left(q^{(k)}\right)
\text{.}
\label{eq:(6)}
\end{equation}

\item[(f)]
For each $k\in \{1,...,L\}$ and $L\ge 1$,
\begin{equation}
\mathcal{G}\left(q^{(0)},q^{(k)}\right)
\le \rho_{1}+\cdots+\rho_{k}
< (Q-1)\cdot\rho_{k}
\text{.}
\label{eq:q^0-q^k}
\end{equation}

\end{enumerate}

\end{proposition}

\proof

For a given
$q=\sum\limits_{i=1}^{I} \ell_i [\![q_i]\!]\in\mathbf{Q}_{Q}(\mathbb{R}^n)$,
denote by
$q^{(0)}=q$ and
$q^{(0)}
=\sum\limits_{i=1}^{I_0} \ell_i^0 [\![q_i^0]\!]$.
Without loss of generality, we may assume that
$\text{card}\left(\text{spt}(q)\right)\ge2$.
Then we follow the so-called
\emph{standard modification procedure for members of
$\mathbf{Q}_{Q}(\mathbb{R}^n)$}
in {\cite[2.9]{Alm00}} to find $q^{(1)}$
by choosing sufficiently large constant $K$ in {\cite[2.9]{Alm00}}
(see (\ref{eq:K-def}) below).
We successively apply this procedure to obtain $q^{(k+1)}$
from $q^{(k)}$ until
$\text{card}\left(\text{spt}(q^{L})\right)=1$ for some $L\in\mathbb{Z}_{+}$.
Such positive integer $L\le Q-1$ fulfilling
$\text{card}\left(\text{spt}(q^{(L)})\right)=1$
exists, because we will show that 
$\text{card}\left(\text{spt}(q^{(k)})\right)
\gvertneqq \text{card}\left(\text{spt}(q^{(k+1)})\right)$
in the procedure.

Note that the assertion in Proposition is obviously fulfilled as $k=0$.
Thus, by induction argument,
we suppose that the assertion is true for
$q^{(0)},...,q^{(k)}\in\mathbf{Q}_{Q}(\mathbb{R}^n)$,
where $k\ge 1$ and $\text{card}\left(\text{spt}(q^{(k)})\right)\ge2$
(otherwise, the proof is finished).
Now, we show how to find $q^{(k+1)}$ from $q^{(k)}$. 
We first define some numbers for the construction of
the sequence in (\ref{eq:(rho-sigma)-sequence}) as follows.
Let
\begin{equation}
K=K(n,Q)=\frac{20 Q}{\sin\theta_{0}(n,Q)} 
\label{eq:K-def}
\end{equation}
and
\begin{equation}
s_{0}=\sigma_{k} 
\label{eq:s_(0)=sigma_(k)}
\end{equation}
where $\sigma_{k}$ is given in (\ref{eq:sigma_(k)-def}).
The choice of $K$ in (\ref{eq:K-def}) assures that the
quotient $\sigma_{k}/\rho_{k}$ is sufficient large,
see (\ref{eq:Lower_Bdd-sigma_(k)/rho_(k)}).
Besides,
for a constructed $q^{(0)},...,q^{(k)}$,
a sufficiently large $K$ and a properly chosen $\sigma_{k}$ 
assure that 
$\text{card}\left(\text{spt}(q^{(k)})\right)
\gvertneqq\text{card}\left(\text{spt}(q^{(k+1)})\right)$,
in the following procedure.
Let
\begin{equation}
\left\{
\begin{array}{l}
t_{1}=0 \text{,}
\\
d_{\varkappa}=(Q-1)t_{\varkappa} \text{, for each }
\varkappa=1,2,...\text{,}
\\
s_{\varkappa}=(Q-1)d_{\varkappa}+s_{\varkappa-1} \text{, for each }
\varkappa=1,2,...\text{,}
\\
t_{\varkappa+1}=2K\cdot s_{\varkappa} \text{, for each }
\varkappa=1,2,...\text{.}
\end{array}
\right.
\label{def:The-sequence-(modification_procedure)-1}
\end{equation}
By a direct computation from
(\ref{def:The-sequence-(modification_procedure)-1})
and restricting $\varkappa\in\{1,...,Q\}$,
we could derive
\begin{equation}
\left\{
\begin{array}{l}
s_{0}=s_{1}  \text{,}
\\
s_{\varkappa}=[1+2K(Q-1)^2]^{\varkappa-1} s_{0}
\le [1+2K(Q-1)^2]^{Q-1} s_{0}
\text{,}
\\
d_{\varkappa}=2K(Q-1)[1+2K(Q-1)^2]^{\varkappa-2} \cdot s_{0}
\text{.}
\end{array}
\right.
\label{def:The-sequence-(modification_procedure)-2}
\end{equation}
For each fixed $\varkappa\in\{1,...,Q\}$, define a partitioning of
$\text{spt}(q^{(k)})=\{q_1^k,...,q_{I_k}^{k}\}$ into equivalence classes by
saying that $q_{i}^{k}\sim q_{j}^{k}$
if there exists a sequence
$\{q_{i_1}^{k},...,q_{i_A}^{k}\}
\subset\text{spt}(q^{(k)})$
such that
$q_{i}^{k}=q_{i_1}^{k}$, $q_{j}^{k}=q_{i_A}^{k}$ and
$|q_{i_\alpha}^{k}-q_{i_{\alpha+1}}^{k}|\le t_k$
for each $\alpha=1,2,...,A-1$.
Denote by $P(\varkappa,1),...,P(\varkappa,N_{k}(\varkappa))
\subset\{q_{1}^{k},...,q_{I_{k}}^{k}\}$
a list of the distinct equivalence classes,
where $N_{k}(\varkappa)$ represents the number of distinct equivalence classes
at the $\varkappa$-th stage of this partitioning procedure.
It is easy to see that
\begin{equation}
\text{diam}\left(P(\varkappa, i)\right)\le d_{\varkappa}
\text{, } \forall \text{ }\varkappa \text{ and } i
\label{eq:Diam(P)}
\end{equation}
and
$Q\ge N_{k}(1)\ge N_{k}(2) \ge \cdots\ge 1$.
Denote by $\varkappa_0$ the smallest positive integer among the integers
$\varkappa$ so that
$N_{k}(\varkappa)=N_{k}(\varkappa+1)$.
It is clear that $\varkappa_0\le Q$.
Furthermore, let
\begin{equation}
\rho_{k+1}:=s_{\varkappa_{0}}
\label{eq:rho_(k+1)-def}
\end{equation}
and, for each $i\in\{1,...,N_{k}(\varkappa_0)\}$,
choose some
$q_{i}^{k+1}\in P(\varkappa_0,i)$
and let
$\ell_{i}^{k+1}=\text{card}\left(P(\varkappa_0,i)\right)$.

From (\ref{eq:Diam(P)}),
it is easy to check that for each fixed integer $k\ge0$,
\begin{equation}
\begin{array}{l}
\left[\mathcal{G}\left(q^{(k)}, q^{(k+1)}\right)\right]^{2}
=\left[\mathcal{G}\left(
\sum\limits_{i=1}^{I_k} \ell^{k}_i [\![q_i^{k}]\!],
\sum\limits_{j=1}^{I_{k+1}} \ell^{k+1}_j [\![q_j^{k+1}]\!]
\right)\right]^{2}
\\
\le (Q-1)\cdot
\left[
\sup\limits_{i}\left\{\text{diam}\left(P(\varkappa_{0},i)\right)\right\}
\right]^{2}
\le  (Q-1)\cdot d_{\varkappa_{0}}^{2}
\text{.}
\end{array}
\label{def:The-sequence-(modification_procedure)-3}
\end{equation}
Moreover, for any $z\in\mathbf{Q}_{Q}(\mathbb{R}^n)$ satisfying
$\mathcal{G}\left(z, q^{(k)}\right) \le s_{0}$,
we may apply the triangle inequality of $\mathcal{G}(\cdot,\cdot)$
in $\mathbf{Q}_{Q}(\mathbb{R}^n)$,
and (\ref{def:The-sequence-(modification_procedure)-2}),
(\ref{def:The-sequence-(modification_procedure)-3})
to derive
\begin{equation}
\begin{array}{l}
\mathcal{G}\left(z, q^{(k+1)}\right)
\le \mathcal{G}\left(z, q^{(k)}\right)
+ \mathcal{G}\left(q^{(k)}, q^{(k+1)}\right)
\\
\le s_{0}+ (Q-1)^{1/2}\cdot d_{\varkappa_{0}}
\\
= \left(1+2K(Q-1)^{3/2}[1+2K(Q-1)^2]^{\varkappa_{0}-2}\right)
\cdot s_{0}
\\
\le
[1+2K(Q-1)^2]^{\varkappa_{0}-1} \cdot s_{0}
=s_{\varkappa_{0}}
\end{array}
\label{def:The-sequence-(modification_procedure)-4}
\end{equation}
where the last equality comes from
(\ref{def:The-sequence-(modification_procedure)-2}).
From (\ref{eq:rho_(k+1)-def}) and (\ref{eq:s_(0)=sigma_(k)}), 
we have 
\begin{equation*}
\mathcal{G}\left(z, q^{(k+1)}\right) \le \rho_{k+1}
%\label{def:The-sequence-(modification_procedure)-5}
\end{equation*}
if $z$ satisfies $\mathcal{G}\left(z, q^{(k)}\right) \le \sigma_{k}$.
Note that, due to the choice of $K$ in (\ref{eq:K-def})
and by letting $s_0=\sigma_{k}$ in (\ref{eq:s_(0)=sigma_(k)}),
we derive from (\ref{def:The-sequence-(modification_procedure)-1}),
\begin{equation*}
t_2=10Q\cdot\inf\limits_{i\ne j}\left\{|q_i^k-q_j^k|\right\}
>\inf\limits_{i\ne j}\left\{|q_i^k-q_j^k|\right\}
\text{.}
\end{equation*}
This implies that $\varkappa_{0}\ge 2$
and therefore
\begin{equation}
\text{card}\left(\text{spt}(q^{(k)})\right)
\gvertneqq
\text{card}\left(\text{spt}(q^{(k+1)})\right)
\text{.}
\label{eq:cardinality-strictly-decreasing}
\end{equation}
Thus, from (\ref{eq:rho_(k+1)-def}) and (\ref{eq:sigma_(k)-def}),
we have
$\sigma_{k}\lvertneqq\rho_{k+1}$.
Besides,
from  (\ref{eq:rho_(k+1)-def}), (\ref{eq:s_(0)=sigma_(k)})
and (\ref{def:The-sequence-(modification_procedure)-2}),
we have
\[
\frac{\rho_{k+1}}{\sigma_{k}}
=[1+2K(Q-1)^2]^{\varkappa_{0}-1}
\le [1+2K(Q-1)^2]^{Q-1}
\text{.}
\]
Now we let
\begin{equation*}
C_{0}(n,Q) :=
[1+2K(Q-1)^2]^{Q-1}
\text{.}
%\label{eq:C_{0}-def}
\end{equation*}
On the other hand, the choice of $\varkappa_{0}$
implies that,
\begin{equation}
|z_{i}-z_{j}| > t_{\varkappa_{0}+1}=2K\cdot s_{\varkappa_{0}}
=2K\cdot\rho_{k+1}
\label{eq:sigma_(k)/rho_(k)-Lower_Bdd-1}
\end{equation}
for any $z_{i}\in P(\varkappa_{0},i)$, $z_{j}\in P(\varkappa_{0},j)$
and any distinct equivalence classes
$P(\varkappa_{0},i)$, $P(\varkappa_{0},j)$.
Thus, 
from (\ref{eq:sigma_(k)-def}) and (\ref{eq:K-def}),
(\ref{eq:sigma_(k)/rho_(k)-Lower_Bdd-1}) implies
\begin{equation*}
\sigma_{k+1} \ge \frac{\sin\theta_0}{4}\cdot t_{\varkappa_{0}+1}
= 10Q\cdot \rho_{k+1}
\text{.}
%\label{eq:sigma_(k)/rho_(k)-Lower_Bdd-2}
\end{equation*}
Notice that,
due to the strictly decreasing of cardinality in
(\ref{eq:cardinality-strictly-decreasing}),
we only follow this procedure in constructing
$q^{(k+1)}$ from $q^{(k)}$ for at most $(Q-1)$ many times
(until $\text{spt}\left( q^{(L)} \right)$ is
consisted of only one point).
Therefore, $L\in\{1,...,Q-1\}$.

The proof of each statement from (a) to (e) follows
from the argument above for all $k\in\{1,...,L\}$ inductively.
The proof of (f) follows from applying (e) and the triangle inequality
of the metric space
$\left(\mathbf{Q}_{Q}(\mathbb{R}^n), \mathcal{G}(\cdot,\cdot)\right)$. 

\endproof

\begin{definition}
[The nested admissible closed balls of $q\in\mathbf{Q}_{Q}(\mathbb{R}^n)$]
\label{def:Nested_Adm_Closed_Ball}
We follow the notations in Proposition \ref{prop:Nested_Adm_Closed_Ball}.
Let
$\left\{\mathbb{B}^{\mathbf{Q}}_{\tau_{k}}(q^{(k)})\right\}_{k=0}^{L}$
fulfill
$\tau_{k}\in [\rho_{k},\sigma_{k}]$ for each $k\in\{0,1,...,L\}$,
as stated in Proposition \ref{prop:Nested_Adm_Closed_Ball}.
Then,
$\left\{\mathbb{B}^{\mathbf{Q}}_{\tau_{k}}(q^{(k)})\right\}_{k=0}^{L}$
is said to be a sequence of nested admissible closed balls of $q$.
\end{definition}

The nested admissible closed balls of any member in
$\mathbf{Q}_{Q}(\mathbb{R}^n)$ will be useful
for establishing the ``global" monotonicity formula in the proof of
Lemma \ref{lem:prethm}.

\begin{lemma}[Key Lemma]
\label{lem:prethm}
Assume
$\Omega \subset \mathbb{R}^{2}$ is an open set,
$\mathbb{U}_{r}(w) \subset\subset \Omega$
is an open ball of radius $r>0$ with the center $w$
and
$f\in \mathcal{Y}_{2}(\Omega,\mathbf{Q}_{Q}(\mathbb{R}^n))$
is weakly stationary-harmonic.
Let $h:\Omega\rightarrow\mathbb{R}^{2}$
be the harmonic function induced from
the Hopf-differential of $\boldsymbol\xi_{0}\circ f$
as constructed in Proposition \ref{prop:wc}.
Suppose
$A\subset \Omega$ is a set of Lebesgue points of
$|\nabla (\boldsymbol\xi_{0}\circ f)|^2$
and $w^{\ast}\in A \cap \mathbb{U}_{r}(w)$
fulfills
$\underset{x\in\partial\mathbb{U}_r(w)}{\inf}\mathcal{G}(f(x),f(w^{\ast}))>0$.
Then,
\begin{equation}
\underset{x\in \partial \mathbb{U}_r (w)}
{\inf}
\mathcal{G}(f(x),f(w^{\ast}))
\le
\sqrt{
\frac{Dir \left( f; \mathbb{U}_{r} (w) \right)+ Dir \left(h;\mathbb{U}_{r}(w) \right)}
{2\pi \cdot \delta (n,Q)}
}
\label{eq:Bdd_Osc(f)}
\end{equation}
for some constant $\delta(n,Q)>0$.
\end{lemma}

\begin{remark}
\label{rk:good_point_G}
The set $A\subset\Omega$,
defined as the Lebesgue point of
$|\nabla (\boldsymbol\xi_{0}\circ f)|^2$ in Lemma \ref{lem:prethm},
is independent of the choice of $\boldsymbol\xi_{0}$
(since $|\nabla (\boldsymbol\xi_{0}\circ f)|$
is invariant with respect to the choice of $\boldsymbol\xi_{0}$).
Besides,
$\mathcal{L}^{2}(\Omega \setminus A)=0$,
because $\boldsymbol\xi_{0}\circ f \in W^{1, 2}(\Omega,\mathbb{R}^{N})$.
Recall from Proposition \ref{prop:wc} that
$G=(\boldsymbol\xi_{0}\circ f, h)$ and
\begin{equation}
\lvert \nabla G\rvert^{2}
=\lvert \nabla (\boldsymbol\xi_{0}\circ f) \rvert^{2}
+\lvert \nabla h\rvert^{2}
=\lvert \nabla (\boldsymbol\xi_{0}\circ f) \rvert^{2}
+\frac{|\varphi|^2}{8}+2
\ge 2
\text{.}
\label{eq:nabla_G>0}
\end{equation}
From (\ref{eq:nabla_G>0})
and the definition of ``good" points in Definition \ref{def:good_point_sets_0},
$A\subset\Omega$ is a set of ``good" points of $G$.
\end{remark}

For our convenience, let
$F=\boldsymbol\xi_{0}\circ f \in \mathbb{R}^{nQ}$,
$G=(F, h)\in\mathbb{R}^{nQ+2}$, and
\begin{equation}
d^\ast(x):=\sqrt{\mathcal{G}(f(w^{\ast}), f(x))^2 + |h(w^{\ast})-h(x)|^2}
\text{.}
\label{eq:d^ast}
\end{equation}
Therefore, (\ref{eq:Bdd_Osc(f)}) is equivalent to
\begin{equation}
\underset{x\in \partial \mathbb{U}_r (w)}
{\inf}
d^\ast(x)
\le
\sqrt{
\frac{Dir \left( G; \mathbb{U}_{r} (w) \right)}
{2\pi \cdot \delta (n,Q)}
}
\text{.}
\label{eq:Bdd_Osc(f)_1}
\end{equation}

\proof

Let
$\mathbb{U}_{R_0}(w)\subset\Omega$,
where $w\in\Omega$ and $R_{0}>0$.
Without loss of generality, assume that $f$ is not identically the constant
$f (w^{\ast})$ on $\partial B_{r} (w)$.
Let
\begin{equation*}
\tau_\ast :=
\underset{x\in\partial\mathbb{U}_{r} (w)}{\inf}
\mathcal{G}(f(w^{\ast}), f(x)) >0
\text{.}
\end{equation*}
We follow the notations in Proposition \ref{prop:Nested_Adm_Closed_Ball}
and let $f (w^{\ast})=q^{(0)}$,
where $w^{\ast}\in A\cap\mathbb{U}_{r}(w)$.
From Proposition \ref{prop:Adm_Ball_Proj},
we may also choose a suitable coordinate bases of $\mathbb{R}^{n}$ 
for the construction of admissible closed balls of
$w^{\ast}$.
Since the Dirichlet integral is invariant under the change of
Cartesian coordinates of $\mathbb{R}^{n}$,
for our convenience, we still denote by $f$
the multiple-valued function
after composed with the change of coordinates of $\mathbb{R}^{n}$.
Then,
from Proposition \ref{prop:Nested_Adm_Closed_Ball},
there correspond two sequences of non-negative numbers
$\{\rho_{k}\}_{k=0}^{L}$,
$\{\sigma_{k}\}_{k=0}^{L}$ and a sequence of members
$\{q^{(k)}\}_{k=0}^{L}\subset \mathbf{Q}_{Q}(\mathbb{R}^{n})$.
For a given positive number $\tau_{\ast}$, let
$k_{0} \in \{0,1,...,L\}$ be the integer fulfilling one of the following conditions,
\begin{equation}
\left\{
\begin{array}{l}
10 Q\cdot\rho_{k_{0}}\lvertneqq\tau_{\ast}\le 10 Q\cdot\rho_{k_{0}+1}
\text{, } k_{0}\in\{0,1,...,L-1\}\text{, } 
\\
10 Q\cdot\rho_{L}\lvertneqq\tau_{\ast}< +\infty
\text{, } k_{0}=L \text{.}
\end{array}
\right.
\label{eq:def-k_0}
\end{equation}
Let
\begin{equation}
d^{\ast}_{k}(x):=\sqrt{\mathcal{G}(q^{(k)},f(x))^2 + |h(w^{\ast})-h(x)|^2}
\label{eq:d^ast_k}
\end{equation}
and
\begin{equation}
\Omega^{\ast}_{k}(\rho):=
\left\{
x\in\Omega:
d^\ast_{k}(x)<\rho
\right\}
\label{eq:Omega_k}
\end{equation}
for each $k\in\{0,1,...,L\}$.
For any $x$ satisfying one of the following conditions
\begin{equation}
\left\{
\begin{array}{l}
d^\ast_{k}(x)
\le
\sigma_{k}
\text{, if }
k\in\{0,...,k_{0}-1\}
\text{,}
\\
d^\ast_{k_0}(x)
\le
\frac{2}{5}\min\{\tau_{\ast},\sigma_{k_0}\}
\text{,}
\end{array}
\right.
\label{eq:Upper_Bdd-d^ast_k}
\end{equation}
one could verify that
\begin{equation}
\begin{array}[t]{l}
d^\ast_{0}(x)
:= \sqrt{\mathcal{G}(q^{(0)},f(x))^2 + |h(w^{\ast})-h(x)|^2}
\le \mathcal{G}\left(q^{(0)}, f(x)\right) + |h(w^{\ast})-h(x)|
\\ \\
\le
\mathcal{G}\left(q^{(0)}, q^{(k)}\right) + \mathcal{G}\left(q^{(k)}, f(x)\right)
+ |h(x)-h(w^{\ast})|
%\\ \\
\le
\mathcal{G}\left(q^{(0)}, q^{(k)}\right) + \sqrt{2} \cdot d^\ast_{k}(x)
\text{.}
\end{array}
\label{eq:cpt_spt-0}
\end{equation}
As $k\in\{0,...,k_{0}-1\}$,
we apply (\ref{eq:Upper_Bdd-d^ast_k}),
(\ref{eq:(rho-sigma)-sequence}), (\ref{eq:q^0-q^k}) and (\ref{eq:def-k_0})
to derive
\begin{equation*}
\begin{array}[t]{l}
\text{R.H.S. of (\ref{eq:cpt_spt-0})}
\lvertneqq
\mathcal{G}\left(q^{(0)},q^{(k)}\right)
+ \sqrt{2}\cdot\sigma_{k}
\\ \\
<(Q-1)\cdot\rho_{k}+ \sqrt{2}\cdot \sigma_{k}
<(1+\sqrt{2})\cdot \sigma_{k}
< (1+\sqrt{2})\cdot \rho_{k_0} < \tau_{\ast}
\text{.}
\end{array}
%\label{eq:cpt_spt-1}
\end{equation*}
As $k=k_0$, we apply (\ref{eq:Upper_Bdd-d^ast_k}) to derive
\begin{equation*}
\begin{array}[t]{l}
\text{R.H.S. of (\ref{eq:cpt_spt-0})}
\lvertneqq
\mathcal{G}\left(q^{(0)},q^{(k_0)}\right)
+ \frac{2\sqrt{2}}{5}\cdot\min\{\tau_{\ast},\sigma_{k_0}\}
\\ \\
=
\min\{\tau_{\ast},\sigma_{k_0}\}
-
\left(
\frac{5-2\sqrt{2}}{5}\cdot\min\{\tau_{\ast},\sigma_{k_0}\}
- \mathcal{G}\left(q^{(0)},q^{(k_0)}\right)
\right)
\\ \\
\lvertneqq
\min\{\tau_{\ast},\sigma_{k_0}\}
\text{,}
\end{array}
%\label{eq:cpt_spt-1}
\end{equation*}
where the last inequality comes from applying
(\ref{eq:def-k_0}), (\ref{eq:Lower_Bdd-sigma_(k)/rho_(k)}) and (\ref{eq:q^0-q^k}).
Therefore, 
\begin{equation}
\left\{
\begin{array}{l}
%\\ \\
\rho \in \left(0, \sigma_{k}\right)
\text{, as } k=0,...,k_{0}-1
\text{,}
\\ \\
\rho \in \left(0, \frac{2}{5}\min\{\tau_{\ast},\sigma_{k_{0}}\}\right)
\text{, as } k=k_{0}
\text{.}
\end{array}
\right.
\Longrightarrow
\Omega^{\ast}_{k}(\rho)\subset\subset\mathbb{U}_{r}(w)
\text{.}
\label{eq:Omega^ast_(k)subsetsubsetU_r}
\end{equation}

In the rest of this article, for a given $k\in\{0,...,k_0\}$,
we always let $\rho$ fulfill the condition in
(\ref{eq:Omega^ast_(k)subsetsubsetU_r}).

\bigskip

\textbf{Step 1: Constructing ``admissible" range-variations.}

Define the smooth vector field
$\Gamma_{k}:\mathbb{R}^{n}\rightarrow\mathbb{R}^{n}$
by
\begin{equation*}
\Gamma_{k}(y):=
\chi(|q_{i}^{k}-y|)\cdot\left(q_{i}^{k}-y\right) 
\end{equation*}
where
$\chi:\mathbb{R}\rightarrow[0,1]$ is a
smooth function satisfying
$\chi(s)=1$ as $s\le \frac{2}{5}\sigma_{k}$,
$\chi(s)=0$ as $s\ge  \frac{3}{5}\sigma_{k}$.
Thus, $\text{Lip}(\Gamma_{k})\le \frac{5}{\sigma_{k}}$. 
The $\mathbf{Q}_{Q}(\mathbb{R}^n)$-valued function,
induced from the smooth vector field
$\Gamma_{k}:\mathbb{R}^{n}\rightarrow\mathbb{R}^{n}$,
can be written as
\begin{equation}
\left(\Gamma_{k}\right)_{\#}(f(x))
:= \sum\limits_{i=1}^{Q} [\![ \Gamma_{k}\circ f_{i}(x)
]\!]
\text{.}
\label{eq:(Gamma_k)_(star)}
\end{equation}
Note, the definition of $\Omega_{k}^{\ast}(\rho)$ in (\ref{eq:Omega_k})
and the definition of $d_{k}^{\ast}$ in (\ref{eq:d^ast_k})
give us the condition
\begin{equation}
\rho \le \frac{2}{5}\sigma_{k}
\Longrightarrow
f \left( \Omega^\ast_{k}(\rho) \right)
\subset \mathbb{B}^{\mathbf{Q}}_{\frac{2}{5}\sigma_{k}}(q^{(k)})
\text{.}
\label{eq:f(Omega)-subset-B^Q}
\end{equation}
Let
$\lambda \in C^{\infty}(\mathbb{R},[0,1])$
satisfy 
$\lambda^{\prime}(s) \ge 0$
and
\begin{equation}
\lambda (s)=
\left\{
\begin{array}{l}
0 \text{ , if } s \le 0 \text{,}
\\
1 \text{ , if } s \ge \varepsilon \text{,}
\end{array}
\right. 
\label{eq:lambda-def}
\end{equation}
for some $\varepsilon>0$.
From (\ref{eq:Range_Variations}),
we let the range-variation of $f$ be
\begin{equation}
f^{t}(x):= \sum\limits_{i=1}^{Q} [\![ f_{i}(x)
+ t \cdot \lambda(\rho-d^\ast_{k}(x)) \cdot \Gamma_{k}(f_{i}(x)) ]\!] 
\label{eq:The_Range_Variation-f}
\end{equation}
where
$t\in(-\epsilon,\epsilon)$
and $\epsilon>0$ is a sufficiently small number.
It remains to show that the range-variation of $f$
in (\ref{eq:The_Range_Variation-f}) is admissible.
In other words, for each fixed $t$, we should prove that
$\{f^t\}$ belongs to the class of Sobolev space
$\mathcal{Y}_{2}(\Omega,\mathbf{Q}_{Q}(\mathbb{R}^n))$
and $\{f^t\}$ could be approximated by a sequence of multiple-valued functions
generated from a smooth perturbation of $f$
(as defined in Definition \ref{def:WH}).
In the following, we want to show the local fine-property of $\Lambda_{k}$ and $f^t$.

We first prove that
$\Lambda_{k}:\Omega\rightarrow\mathbb{R}$
belongs to the class of Sobolev space
$W^{1, 2}_{0}(\Omega)\cap L^{\infty}(\Omega)$
by the difference quotient method.
Observe that
\begin{equation*}
\begin{array}[t]{l}
\frac{|\Lambda_{k}(x_2)-\Lambda_{k}(x_1)|}{|x_2-x_1|}
=\frac{|\Lambda_{k}(\rho-d^\ast_{k}(x_2))- \Lambda_{k}(\rho-d^\ast_{k}(x_1))|}
{|\left(\rho-d^\ast_{k}(x_2)\right)-\left(\rho-d^\ast_{k}(x_1)\right)|}
\cdot
\frac{|d^\ast_{k}(x_2)-d^\ast_{k}(x_1)|}{|x_2-x_1|}
\\ \\
\le
\text{Lip}(\lambda)\cdot
\frac{|
(\mathcal{G}(q^{(k)},f(x_2))^2-\mathcal{G}(q^{(k)},f(x_1))^2)|
+||h(w^{\ast})-h(x_2)|^2 - |h(w^{\ast})-h(x_1)|^2|}
{|x_2-x_1|\cdot |d^\ast_{k}(x_2)+d^\ast_{k}(x_1)|}
\\ \\
\le
\text{Lip}(\lambda)\cdot
\left(
\frac{|\mathcal{G}(q^{(k)}, f(x_2))-\mathcal{G}(q^{(k)},f(x_1))|}{|x_2-x_1|}
+
\frac{\lvert|h(w^{\ast})-h(x_2)| - |h(w^{\ast})-h(x_1)|\rvert}
{|x_2-x_1|}
\right)
\\ \\
\le
\text{Lip}(\lambda)\cdot
\left(
\frac{\mathcal{G}(f(x_2),f(x_1))}{|x_2-x_1|}+\frac{|h(x_2)-h(x_1)|}{|x_2-x_1|}
\right)
\\ \\
\le
\text{Lip}(\lambda)\cdot
\left(
\text{Lip}(\boldsymbol\xi^{-1})\cdot
\frac{|\boldsymbol\xi\circ f(x_2)-\boldsymbol\xi\circ f(x_1)|}{|x_2-x_1|}
+\frac{|h(x_2)-h(x_1)|}{|x_2-x_1|}
\right)
\end{array}
\end{equation*}
where the third inequality comes from applying triangle inequality of
the metric $\mathcal{G}(\cdot,\cdot)$.
Because $\boldsymbol\xi\circ f\in W^{1, 2}(\Omega)$
and $h$ is a smooth harmonic function,
we derive the uniform bound of the difference quotient
$\frac{\left|\Lambda_{k}(x_2)-\Lambda_{k}(x_1)\right|}{|x_2-x_1|}$
in $L^{2}(\Omega)$-topology for any distinct
$x_2, x_1$ in $\Omega$.
Thus,
$\Lambda_{k}\in W^{1, 2}(\Omega)$.
The restriction of $\rho$ in (\ref{eq:Omega^ast_(k)subsetsubsetU_r})
and the condition of $\Lambda_{k}\in[0,1]$
imply that
$\Lambda_{k}\in W^{1, 2}_{0}(\Omega)\cap L^{\infty}(\Omega)$.

Now we want to prove that
$f^t$ belongs to the class of Sobolev space
$\mathcal{Y}_{2}(\Omega, \mathbf{Q}_{Q}(\mathbb{R}^n))$
for each fixed $t$.
Without loss of generality, we only need to show the case when we let
$\Omega=\mathbb{U}_{r}(w)$ because $f^{t}(x)=f(x)$
for all $t$ and $x\in\Omega\setminus\mathbb{U}_{r}(w)$.
We re-write the equation (\ref{eq:The_Range_Variation-f}) as
\begin{equation}
f^{t}(x)= \sum\limits_{i=1}^{Q} [\![ f^{t}_{i}(x) ]\!]
=\sum\limits_{i=1}^{Q} [\![ f_{i}(x) + t \cdot \Lambda_{k}(x)
\cdot \Gamma_{k}(f_{i}(x)) ]\!]
\text{.}
\label{eq:The_Range_Variation-f-2}
\end{equation}
For fixed $x_1, x_2 \in \mathbb{U}_{r}(w)$, denote by
$\sigma_{x_1, x_2}:\{1,2,...,Q\}\rightarrow\{1,2,...,Q\}$
the permutation fulfilling
\begin{equation}
\begin{array}[t]{l}
\mathcal{G}\left(f(x_2), f(x_1)\right)
= \sqrt{
\sum\limits_{\ell=1}^{Q}
\left|
f_{\sigma_{x_1, x_2}(\ell)}(x_2) - f_{\ell}(x_1)
\right|^{2}
}
\text{.}
\end{array}
\label{eq:sigma-Q_(dist)}
\end{equation}
Then, for fixed $x_1$, $x_2$, $t$ and $\alpha$,
\begin{equation}
\begin{array}[t]{l}
\left| \xi_{\alpha}\circ f^{t}(x_2) - \xi_{\alpha}\circ f^{t}(x_1) \right|^{2}
= \left[\mathcal{G}
\left(
(\Pi_{\alpha})_{\#}(f^{t}(x_2))
,
(\Pi_{\alpha})_{\#}(f^{t}(x_1))
\right)\right]^{2}
\\ \\
= \inf\limits_{\sigma\in \mathcal{P}_{Q}}
\left(
\sum\limits_{\ell=1}^{Q}
\left|
\Pi_{\alpha}\left( f^{t}_{\sigma(\ell)}(x_2) \right)
- \Pi_{\alpha}\left( f^{t}_{\ell}(x_1) \right)
\right|^{2}
\right)
\\ \\
\le
\sum\limits_{\ell=1}^{Q}
\left|
\Pi_{\alpha}\left( f^{t}_{\sigma_{x_1, x_2}(\ell)}(x_2) \right)
- \Pi_{\alpha}\left( f^{t}_{\ell}(x_1) \right)
\right|^{2}
=
\sum\limits_{\ell=1}^{Q}
\left|
\Pi_{\alpha}
\left(
f^{t}_{\sigma_{x_1, x_2}(\ell)}(x_2)
- f^{t}_{\ell}(x_1)
\right)
\right|^{2} 
\end{array}
\label{eq:xi_(alpha)-(f(x_2)-f(x_1))-1}
\end{equation}
where the first equality comes from (\ref{eq:Q(R)-topology}) and
$\mathcal{P}_{Q}$ denotes the permutation group of $\{1,2,...,Q\}$.
By applying the inequality,
$(a+b+c)^2\le 3(a^2+b^2+c^2)$, we have
\begin{equation}
\begin{array}[t]{l}
\sum\limits_{\ell=1}^{Q}
\left|
\Pi_{\alpha}
\left(
f^{t}_{\sigma_{x_1, x_2}(\ell)}(x_2)
- f^{t}_{\ell}(x_1)
\right)
\right|^{2}
\\ \\
\le
\sum\limits_{\ell=1}^{Q}
3 \cdot
\left|
\Pi_{\alpha}
\left(
f_{\sigma_{x_1, x_2}(\ell)}(x_2) - f_{\ell}(x_1)
\right)
\right|^{2}
\\ \\
+
\sum\limits_{\ell=1}^{Q}
3t^2 \cdot \left( \Lambda_{k}(x_2) - \Lambda_{k}(x_1) \right)^2 \cdot
\left|
\Pi_{\alpha}
\left(
\Gamma_{k}\circ f_{\sigma_{x_1, x_2}(\ell)}(x_2)
\right)
\right|^{2}
\\ \\
+
\sum\limits_{\ell=1}^{Q}
3 t^2 \cdot \left(\Lambda_{k}(x_1)\right)^2 \cdot
\left|
\Pi_{\alpha}
\left(
\Gamma_{k}\circ f_{\sigma_{x_1, x_2}(\ell)}(x_2)
- \Gamma_{k}\circ f_{\ell}(x_1)
\right)
\right|^{2}
\text{.}
\end{array}
\label{eq:xi_(alpha)-(f(x_2)-f(x_1))-2}
\end{equation}
Thus, from (\ref{eq:xi_(alpha)-(f(x_2)-f(x_1))-1}),
(\ref{eq:xi_(alpha)-(f(x_2)-f(x_1))-2}), (\ref{eq:sigma-Q_(dist)})
and taking the sum $\sum\limits_{\alpha=1}^{n}$,
we have
\begin{equation*}
\begin{array}[t]{l}
\left|\boldsymbol\xi_{0}\circ f^{t}(x_2)
-\boldsymbol\xi_{0}\circ f^{t}(x_1)\right|^{2}
\\ \\
\le
\sum\limits_{\ell=1}^{Q}
3 \cdot
\left|
f_{\sigma_{x_1, x_2}(\ell)}(x_2) - f_{\ell}(x_1)
\right|^{2}
\\ \\
+
\sum\limits_{\ell=1}^{Q}
3t^2 \cdot \left( \Lambda_{k}(x_2) - \Lambda_{k}(x_1) \right)^2 \cdot
\left|
\Gamma_{k}\circ f_{\sigma_{x_1, x_2}(\ell)}(x_2)
\right|^{2}
\\ \\
+
\sum\limits_{\ell=1}^{Q}
3 t^2 \cdot \left(\Lambda_{k}(x_1)\right)^2 \cdot
\left|
\Gamma_{k}\circ f_{\sigma_{x_1, x_2}(\ell)}(x_2)
- \Gamma_{k}\circ f_{\ell}(x_1)
\right|^{2}
\\ \\
\le
3 \cdot
\left[
\mathcal{G}\left(f(x_2), f(x_1)\right)
\right]^{2}
\\ \\
+
\sum\limits_{\ell=1}^{Q}
3t^2 \cdot \left( \Lambda_{k}(x_2) - \Lambda_{k}(x_1) \right)^2 \cdot
\left|
\Gamma_{k}\circ f_{\sigma_{x_1, x_2}(\ell)}(x_2)
\right|^{2}
\\ \\
+
3 t^2 \cdot \left(\Lambda_{k}(x_1)\right)^2
\cdot |\text{Lip}(\Gamma_{k})|^2 \cdot
\left[
\mathcal{G}\left(f(x_2), f(x_1)\right)
\right]^{2}
\text{.}
\end{array}
\end{equation*}
Thus,
\begin{equation*}
\begin{array}[t]{l}
\frac{
\left| \boldsymbol\xi_{0}\circ f^{t}(x_2)-\boldsymbol\xi_{0}\circ f^{t}(x_1)\right|^{2}
}
{|x_2-x_1|^2}
\\ \\
\le
\frac{3 \left[\mathcal{G}\left(f(x_2), f(x_1)\right) \right]^2}{|x_2-x_1|^2}
+
\frac{3t^2 \left(\Lambda_{k}(x_2)-\Lambda_{k}(x_1)\right)^2}{|x_2-x_1|^2}
Q \left\|\Gamma_{k}\right\|_{L^\infty}^{2}
+ 3 t^2 [\text{Lip}(\Gamma_{k})]^2
\frac{\left[\mathcal{G}\left(f(x_2), f(x_1)\right) \right]^2}{|x_2-x_1|^2}
\\ \\
\le 3 \left(1+t^2 [\text{Lip}(\Gamma_{k})]^2\right)
[\text{Lip}(\boldsymbol\xi^{-1})]^2
\frac{\left|\boldsymbol\xi\circ f(x_2)-\boldsymbol\xi\circ f(x_1)\right|^2}
{|x_2-x_1|^2}
+
\frac{3t^2 \left(\Lambda_{k}(x_2)-\Lambda_{k}(x_1)\right)^2}{|x_2-x_1|^2}
Q \left\|\Gamma_{k}\right\|_{L^\infty}^{2}
\text{,}
\end{array}
%\label{eq:xi_(alpha)-(f(x_2)-f(x_1))-L^_Bdd}
\end{equation*}
where $\left\|\Gamma_{k}\right\|_{L^\infty} \le 2 \sigma_{k}/5$
and
$\text{Lip}(\Gamma_{k}) \le 5/\sigma_{k}$.
Since $k$ is fixed and both $\boldsymbol\xi\circ f$ and $\Lambda_{k}$
belong to the class of Sobolev spaces
$W^{1, 2}$, we conclude that
$\frac{\left| \xi_{\alpha}\circ f^{t}(x_2)-\xi_{\alpha}\circ f^{t}(x_1)\right|}
{|x_2-x_1|}$
is uniformly bounded in $L^2$ for all distinct $x_1$, $x_2$ in
$\mathbb{U}_{r}(w)$ and fixed $t$, $\alpha\in\{1,...,n\}$.
By the difference quotient method, we conclude that
$f^t$ belongs to the Sobolev space
$\mathcal{Y}_{2}\left(\mathbb{U}_{r}(w),
\mathbf{Q}_{Q}(\mathbb{R}^{n})\right)$ for each fixed $t$.

We may also follow the same argument above to show that
the induced multiple-valued function,
$\left(\Lambda_{k}(x)\cdot \left(\Gamma_{k}\right)_{\#}(f(x))\right)
=\sum\limits_{i=1}^{Q} [\![ \Lambda_{k}(x)\cdot\Gamma_{k}(f_{i}(x)) ]\!] $,
belongs to the Sobolev space
$\mathcal{Y}_{2}\left(\Omega, \mathbf{Q}_{Q}(\mathbb{R}^{n})\right)$
and we write 
\[
\text{ap } A_x \left(\Lambda_{k}(x)\cdot \left(\Gamma_{k}\right)_{\#}(f(x))\right)
=\sum\limits_{i=1}^{Q} [\![ \text{ap }D_{x}
\left(\Lambda_{k}(x)\cdot\Gamma_{k}(f_{i}(x))\right) ]\!]
\text{ for } \mathcal{L}^{2} \text{ a.e. } x \in\Omega
\text{.}
\]

Below we show that
$\frac{d}{dt}_{\vert t=0}Dir (f^{t}; \Omega)=0$
by comparing
$\frac{d}{dt}_{\vert t=0}Dir (f^t; \Omega)$
with
$\frac{d}{dt}_{\vert t=0}Dir (\overset{\sim}{f^t}; \Omega)$,
where $\overset{\sim}{f^t}$ is a smooth perturbation of $f$
(see (\ref{eq:Lusin-1}), (\ref{eq:Lusin-2})
and (\ref{eq:The_Range_Variation-f_(delta)}) for details). 
Note, $f$, $f^t$ and $\Lambda_{k}(x)\cdot \left(\Gamma_{k}\right)_{\#}(f(x))$
all belong to the Sobolev space
$\mathcal{Y}_{2}(\Omega,\mathbf{Q}_{Q}(\mathbb{R}^n))$,
i.e.,
$\boldsymbol\xi\circ f$, $\boldsymbol\xi\circ f^t$ and
$\boldsymbol\xi \left(\Lambda_{k}(x)\cdot \left(\Gamma_{k}\right)_{\#}(f(x))\right)$ all belong to
$W^{1, 2} (\Omega)$.
By the fine-property of functions in Sobolev spaces
(e.g., see {\cite[6.1.3]{EG92}}),
the approximate derivatives of
$\boldsymbol\xi\circ f$, $\boldsymbol\xi\circ f^t$ and
$\boldsymbol\xi \left(\Lambda_{k}(x)\cdot \left(\Gamma_{k}\right)_{\#}(f(x))\right)$ all
exist $\mathcal{L}^2$ a.e. in $\Omega$ (e.g., see {\cite[6.1.3]{EG92}}).
Recall from the description on the
\texttt{strongly approximately affinely approximable}
multiple-valued functions in the Preliminaries of this article
(or see {\cite[1.4 (3)]{Alm00}},
if $\boldsymbol\xi\circ f$ is approximate differentiable at $x$
and $f_{i}(x)=f_{j}(x)$, then
\begin{equation*}
\left\{
\begin{array}{l}
\text{ap } D_{x} f_{i}(x) = \text{ap } D_{x} f_{j}(x) 
\text{,}
\\
\text{ap } D_{x}\left(\Lambda_{k}(x) \cdot \Gamma_{k}\circ f_{i}(x) \right)
=\text{ap } D_{x}\left(\Lambda_{k}(x) \cdot \Gamma_{k}\circ f_{j}(x) \right) 
\text{.}
\end{array}
\right.
\end{equation*}
Thus, if $\boldsymbol\xi\circ f$ 
is approximate differentiable at $x$ and 
$
f(x)=\sum\limits_{\kappa=1}^{K} \text{ } \ell_{\kappa}
\cdot [\![ f_{\kappa}(x) ]\!] 
$, 
then
\begin{equation*}
\left\{
\begin{array}{l}
\text{ap }A_{x} f (x)
=
\sum\limits_{\kappa=1}^{K} \text{ } \ell_{\kappa}\cdot
[\![ \text{ap }D_{x} f_{\kappa} (x) ]\!]
\text{, }
\\
\text{ap } A_{x}\left(\Lambda_{k}(x)\cdot
\left(\Gamma_{k}\right)_{\#}(f(x))\right)
=
\sum\limits_{\kappa=1}^{K} \text{ } \ell_{\kappa}\cdot
[\![ \text{ap }D_{x}\left(\Lambda_{k}(x)
\cdot \Gamma_{k}\circ f_{\kappa}(x) \right) ]\!]
\text{.}
\end{array}
\right.
\end{equation*}
Therefore, for $\mathcal{L}^2$ a.e. $x\in\Omega$,
we have
\begin{equation}
\begin{array}{l}
\text{ap } A_{x} f^{t}(x)= \sum\limits_{i=1}^{Q} \text{ }
[\![\text{ap } D_{x} f^{t}_{i}(x) ]\!]
\\
=\sum\limits_{\kappa=1}^{K} \text{ } \ell_{\kappa}\cdot
[\![ \text{ap } D_{x} f_{\kappa}(x)
+ t \cdot\text{ap } D_{x}\left(\Lambda_{k}(x)\cdot\Gamma_{k}\circ f_{\kappa}(x) \right) ]\!]
\\
=\sum\limits_{i=1}^{Q} \text{ }
[\![ \text{ap } D_{x} f_{i}(x)
+ t \cdot\text{ap } D_{x}\left(\Lambda_{k}(x)\cdot\Gamma_{k}\circ f_{i}(x) \right) ]\!]
\text{.}
\end{array}
\label{eq:The_approximate_derivative-f^t}
\end{equation}
%One might still need to put more details in explaining the formula above!!!!!!!!!
%!!!!!!!!!!!!!!!!!!!!!!!!!!!!!!!!!!!!!!!!!!!!!!!!!!!!!!!!!!!!!!!!!!!!!!!!!!!!!!!!!
%!!!!!!!!!!!!!!!!!!!!!!!!!!!!!!!!!!!!!!!!!!!!!!!!!!!!!!!!!!!!!!!!!!!!!!!!!!!!!!!!!
From {\cite[Theorem 2.2]{Alm00}},
we may compute the Dirichlet integral of $f^t$ on an open set
$\Omega$ by
 $\sum\limits_{i=1}^{Q}\| \text{ap } D_{x} f^{t}_{i} \|_{L^2(\Omega)}^2$. 
Thus, from (\ref{eq:The_approximate_derivative-f^t}), we have
\begin{equation*}
\begin{array}[t]{l}
Dir (f^{t}; \Omega)
= \sum\limits_{i=1}^{Q} \text{ }
\int\limits_{\Omega}\left| \text{ap } D_{x} f^{t}_{i}(x) \right|^2 \text{ }dx
\\ \\
=\sum\limits_{i=1}^{Q} \text{ }
\int\limits_{\Omega}\left|
\text{ap } D_{x} f_{i}(x)
+ t \cdot \left[ \text{ap } D_{x} \Lambda_{k}(x) \text{ } \Gamma_{k}( f_{i}(x) )
+ \Lambda_{k}(x) \cdot \text{ap } D_{x} (\Gamma_{k}\circ f_{i})(x) \right]
\right|^2 \text{ }dx
\\ \\
=\sum\limits_{i=1}^{Q} \text{ }
\int\limits_{\Omega}\left|
\text{ap } D_{x} f_{i}(x)
\right|^2 \text{ }dx
\\ \\
+
2t\cdot
\sum\limits_{i=1}^{Q} \text{ }
\int\limits_{\Omega}
\left< \text{ap } D_{x} f_{i}(x)
: \text{ap } D_{x} \Lambda_{k}(x) \text{ } \Gamma_{k}(f_{i}(x))
\right> \text{ }dx
\\ \\
+
2t\cdot
\sum\limits_{i=1}^{Q} \text{ }
\int\limits_{\Omega}
\left< \text{ap } D_{x} f_{i}(x)
:
\Lambda_{k}(x) \cdot \text{ap } D_{y} \Gamma_{k}(f_{i}(x))\text{ }\text{ap } D_{x} f_{i}(x)
\right> \text{ }dx
\\ \\
+
t^{2} \cdot
\sum\limits_{i=1}^{Q} \text{ }
\int\limits_{\Omega}
\left|
\text{ap } D_{x} \Lambda_{k}(x) \text{ } \Gamma_{k}(f_{i}(x))
+ \Lambda_{k}(x) \cdot \text{ap } D_{x} (\Gamma_{k}\circ f_{i})(x)
\right|^2 \text{ }dx
\text{.}
\end{array}
%\label{eq:The_Dir-Integral-f^t}
\end{equation*}
Thus,
\begin{equation}
\begin{array}[t]{l}
\frac{d}{dt}_{\vert t=0}Dir (f^{t};\Omega)
\\ \\
=
2\cdot
\sum\limits_{i=1}^{Q} \text{ }
\int\limits_{\Omega}
\left< \text{ap } D_{x} f_{i}(x)
: \text{ap } D_{x} \Lambda_{k}(x) \text{ } \Gamma_{k}(f_{i}(x))
\right> \text{ }dx
\\ \\
+
2\cdot
\sum\limits_{i=1}^{Q} \text{ }
\int\limits_{\Omega}
\left< \text{ap } D_{x} f_{i}(x)
:
\Lambda_{k}(x) \cdot \text{ap } D_{y} \Gamma_{k}(f_{i}(x))\text{ }\text{ap } D_{x} f_{i}(x)
\right> \text{ }dx
\text{.}
\end{array}
\label{eq:d_t_Dir-Integral-f^t}
\end{equation}
By applying the Lusin-type Theorem on approximating functions in Sobolev spaces
by $C^1$-smooth functions (e.g., see {\cite[Section 6.6 Corollary 2]{EG92}),
for a given $\Lambda_{k}\in W^{1, 2}_{0}(\Omega, \mathbb{R})$
and any $\delta>0$,
there exists a $C^1$-smooth function
$\overset{\sim}{\Lambda}_{k}\in C^{1}_{0}(\Omega)$ such that
\begin{equation}
\|\overset{\sim}{\Lambda}_{k}-\Lambda_{k}\|_{W^{1, 2}(\Omega)}
\le
\delta
\label{eq:Lusin-1}
\end{equation}
and
\begin{equation}
\mathcal{L}^m \left(\{x:\overset{\sim}{\Lambda}_{k}(x)\ne \Lambda_{k} (x)
\text{ or } \nabla_x \overset{\sim}{\Lambda}_{k}(x)\ne \nabla_x \Lambda_{k}(x)\}\right)
\le\delta 
\label{eq:Lusin-2}
\end{equation}
where $\nabla$ denotes the weak differentiation.
Note that if $H\in W^{1,p}_{\text{loc}}$, then
$\nabla_{x}H(x)=\text{ap } D_{x}H(x)$ for
$\mathcal{L}^{2}$ a.e. $x\in\Omega$,
from standard theory of Sobolev spaces
(e.g., see {\cite[p.233 Remark (i)]{EG92}}).
Let
\begin{equation}
\overset{\sim}{f^t}(x):=
\sum\limits_{i=1}^{Q} [\![ f_{i}(x) + t \cdot \overset{\sim}{\Lambda}_{k}(x)
\text{ } \Gamma_{k}(f_{i}(x)) ]\!]
\text{.}
\label{eq:The_Range_Variation-f_(delta)}
\end{equation}
Then, by following the same computation as the one
in deriving (\ref{eq:d_t_Dir-Integral-f^t}),
we have
\begin{equation}
\begin{array}[t]{l}
\frac{d}{dt}_{\vert t=0}Dir (f^{t}; \Omega)
- \frac{d}{dt}_{\vert t=0}Dir (\overset{\sim}{f^t}; \Omega)
\\ \\
=
2\cdot
\sum\limits_{i=1}^{Q} \text{ }
\int\limits_{\Omega}
\left< \text{ap } D_{x} f_{i}(x)
: \left(\nabla_{x} \Lambda_{k}(x)-\nabla_{x} \overset{\sim}{\Lambda}_{k}(x)\right)
\text{ }\Gamma_{k}(f_{i}(x))
\right> \text{ }dx
\\ \\
+
2\cdot
\sum\limits_{i=1}^{Q} \text{ }
\int\limits_{\Omega}
\left< \text{ap } D_{x} f_{i}(x)
:
\left(\Lambda_{k}(x)-\overset{\sim}{\Lambda}_{k}(x)\right) \cdot \nabla_{y}
\Gamma_{k}(f_{i}(x))
\text{ } \text{ap } D_{x} f_{i}(x)
\right> \text{ }dx
\text{.}
\end{array}
\label{eq:Approximating-d_t_Dir-Integral-f^t}
\end{equation}
As $\delta\rightarrow 0^+$,
the first term on the R.H.S. of (\ref{eq:Approximating-d_t_Dir-Integral-f^t})
tends to zero by applying (\ref{eq:Lusin-1}),
the finiteness of both
$\|\Gamma_{k}\|_{L^\infty}$ and $Dir (f; \Omega)$;
meanwhile
the second term on the R.H.S. of (\ref{eq:Approximating-d_t_Dir-Integral-f^t})
also tends to zero by applying (\ref{eq:Lusin-2}),
the finiteness of
$\|\Lambda_{k}\|_{L^\infty}$,
$\|\overset{\sim}{\Lambda}_{k}\|_{L^\infty}$,
$\|\nabla\Gamma_{k}\|_{L^\infty}$
and the Lebesgue-integrability of
$\sum\limits_{i=1}^{Q}\text{ }\left| \text{ap } D_{x} f_{i}\right|^2$. 
Thus, we conclude that
\begin{equation*}
\left| \frac{d}{dt}_{\vert t=0}Dir (f^{t}; \Omega)
- \frac{d}{dt}_{\vert t=0}Dir (\overset{\sim}{f^t}; \Omega) \right|
\rightarrow 0 
\end{equation*}
as $\delta\rightarrow 0^+$.
Since the vector field
$\psi=\overset{\sim}{\Lambda}_{k}\cdot\Gamma_{k}\in C^{\infty}_{c}
(\Omega\times\mathbb{R}^n,\mathbb{R}^n)$ for any $\delta>0$,
from Definition \ref{def:WH} and the assumption of $f$ being weakly harmonic,
we have $\frac{d}{dt}_{\vert t=0}Dir (\overset{\sim}{f^t}; \Omega)=0$.
Therefore, we conclude that
$\frac{d}{dt}_{\vert t=0}Dir (f^t; \Omega)=0$ must hold.

\textbf{Step 2: Deriving the ``local" monotonicity formula.}

%Below needs to be modified!
For simplicity of notations, let
\[
G:=(F, h):=(\boldsymbol\xi_{0}\circ f, h)
\]
and
\[
G^{\ast}_{k}:=\left(\boldsymbol\xi_{0} (q^{(k)}), h(w^\ast)\right)
\text{, }
\forall \text{ } k\in\{0,...,L\}
\text{.}
\]
Notice that, from (\ref{eq:|xi(p)-xi(q)|=G(p,q)}), we have
\begin{equation}
d^{\ast}_{k}(x)=|G^{\ast}_{k}-G(x)|
\text{, } \forall \text{ } x \in \Omega^\ast_{k}\left(\frac{2}{5}\sigma_{k}\right)
\text{.}
\label{eq:d^ast_k=|G*-G|}
\end{equation}
Besides,
due to the choice of cut-off function
$\lambda$ and the definition of
$\Omega_{k}^{\ast}(\rho)$ in (\ref{eq:Omega_k}),
we have
\begin{equation}
\lambda(\rho-d^\ast_{k}(x))=0\text{, }\forall\text{ }x\notin\Omega_{k}^{\ast}(\rho)
\text{.}
\label{eq:Set_of_lambda=0}
\end{equation}
Since both $\lambda$ and $\Gamma_{k}$
in (\ref{eq:The_Range_Variation-f}) are bounded,
we have
\begin{equation*}
f^t(x) \in
\mathbb{B}^{\mathbf{Q}}_{\frac{2}{5}\sigma_{k}}(q^{(k)})
\end{equation*}
for all $x \in \Omega$, $\rho\le\frac{2}{5}\sigma_{k}$
and sufficiently small $|t|$.
Thus, under this situation, there is a unique way of paring
$f^t(x)$ with $q^{(k)}$ in the sense of (\ref{eq:The_subtraction}). 
By applying (\ref{eq:The_linear_relation}), we obtain the expression
\begin{equation*}
\boldsymbol\xi_{0}\circ f^t(x)
= (1-t\cdot\lambda(\rho-d^{\ast}_{k}(x)))\cdot \boldsymbol\xi_{0}\circ f(x)
+ t\cdot\lambda(\rho-d^{\ast}_{k}(x)))\cdot \boldsymbol\xi_{0}\left(q^{(k)}\right)
\text{.}
%\label{eq:The_interpolation-2}
\end{equation*}
In other words, 
\begin{equation*}
F^{t}(x):=
\boldsymbol\xi_0\circ f^t (x)=
F(x)+ t\cdot\lambda(\rho-d^\ast_{k}(x)) \cdot
\left( \boldsymbol\xi_{0} (q^{(k)}) -F(x) \right) 
\end{equation*}
for all $x\in\Omega$, $\rho\le\frac{2}{5}\sigma_{k}$ and sufficiently small $|t|$.
Furthermore, by letting
\begin{equation*}
h^{t}(x):=
h(x)+ t\cdot\lambda(\rho-d^\ast_{k}(x)) \cdot \left( h(w^\ast)-h(x) \right) 
%\label{eq:The_Range_Variation-h}
\end{equation*}
we obtain the perturbation formula of $G$,
\begin{equation}
G^{t}(x):= (F^{t}(x), h^{t}(x)) =
G(x)+ t\cdot\lambda(\rho-d^\ast_{k}(x))
\cdot \left( G^{\ast}_{k}-G(x) \right)
\text{.}
\label{eq:Perturbation_in_Range}
\end{equation}
Note that in the rest of this paper,
the perturbation formula of $f$ and $G$ will be applied only
as $\rho \in (\rho_{k}, \frac{2}{5}\sigma_{k})$ for $k \in \{0,1,..., k_{0}-1\}$ and
$\rho \in (\rho_{k_{0}}, \frac{2}{5}\min\{\tau_{\ast},\sigma_{k_{0}}\})$.
Due to (\ref{eq:Set_of_lambda=0}) and
(\ref{eq:Omega^ast_(k)subsetsubsetU_r}),
the perturbations in (\ref{eq:Perturbation_in_Range}) leave the boundary
value of $G$ fixed, i.e.,
\begin{equation*}
G^{t}(x)=G(x)
\text{, }\forall \text{ }x\in\partial\mathbb{U}_{r}(w) 
\text{.}
%\label{eq:Bdry_Fixed}
\end{equation*}
Since $f$ is a stationary-harmonic multiple-valued function
and $h$ is a harmonic (single-valued) function,
\[
\frac{d}{dt}_{\lfloor t=0} Dir (F^{t};\Omega)
=0
=\frac{d}{dt}_{\lfloor t=0} Dir (h^{t};\Omega)
\text{.}
\]
Thus,
\begin{equation*}
\frac{d}{dt}_{\lfloor t=0} Dir (G^{t};\Omega)
=\frac{d}{dt}_{\lfloor t=0} Dir (f^{t};\Omega)
+ \frac{d}{dt}_{\lfloor t=0} Dir (h^{t};\Omega)
=0
\text{.}
\end{equation*}
From (\ref{eq:d^ast_k=|G*-G|}), this implies that
\begin{equation}
\begin{array}{l}
0 =
-\frac{1}{2} \frac{d}{dt}_{\lfloor t=0} Dir (G^{t}; \Omega)
\\
=
\underset{
\mathbb{U}_r (w)\cap \Omega^{\ast}_{k}(\rho)}{\int }
\lambda \left( \rho- |G^{\ast}_{k}-G(x)| \right)
\cdot \left| \nabla G (x) \right| ^{2} \text{ } dx
\\
-
\underset{\mathbb{U}_r (w)\cap \Omega^{\ast}_{k}(\rho)}{\int }
\frac{\lambda^{\prime} \left( \rho- |G^{\ast}_{k}-G(x)| \right)}
{|G^{\ast}_{k}-G(x)|}
\cdot 
\sum\limits_{i=1}^{2}
\left\langle \partial_{i} G (x), G^{\ast}_{k}-G (x) \right\rangle^{2}
\text{ } dx
\end{array}
\label{eq:diff_ineq_0}
\end{equation}
where $\partial_{i}:=\frac{\partial}{\partial x_i}$.

Below we show how to apply the weak conformality conditions of $G$ 
to (\ref{eq:diff_ineq_0}). 
Let $a, b : \Omega\subset\mathbb{R}^2 \rightarrow \mathbb{R}^N$
be two differentiable vector-valued functions. 
Observe that 
\begin{equation*}
\sum\limits_{i=1}^{2} \left\langle \partial_{i}a, \frac{b}{|b|}\right\rangle^{2}
=
\sum\limits_{i=1}^{2}
\left\langle
\partial_{i}a, \mathbb{P}_{\{\partial_{1} a,\partial_{2} a \}}\left(\frac{b}{|b|}\right)
\right\rangle^{2}
\leq
\sum\limits_{i=1}^{2} \left\langle \partial_{i} a, T_{b} \right\rangle^{2}
\end{equation*}
where
$\mathbb{P}_{\{\partial_{1} a,\partial_{2} a\}}(v)$ denotes the orthogonal projection
of vector $v$ into the two-dimensional plane spanned by
$\partial_{1} a$ and $\partial_{2} a$
and
$T_b:=\frac{\mathbb{P}_{\{\partial_{1} a,\partial_{2} a\}}\left(\frac{b}{|b|}\right)}
{\left|\mathbb{P}_{\{\partial_{1} a,\partial_{2} a\}}\left(\frac{b}{|b|}\right)\right|}$.
From the conformality, 
$\left| \partial_{1}a\right| =\left|\partial_{2}a\right|$
and
$\left\langle \partial_{1}a, \partial_{2}a\right\rangle =0$,
we have
$T_b=\frac{\partial_{1} a}{|\partial_{1} a|}\cdot\cos\theta
+\frac{\partial_{2} a}{|\partial_{2} a|}\cdot\sin\theta$
for some $\theta$.
Thus, 
\begin{equation*}
\sum\limits_{i=1}^{2}
\left\langle \partial_{i}a, T_{b} \right\rangle^{2}
= \left| \partial_{1}a \right|^2
= \frac{1}{2} \left(\left| \partial_{1}a + \partial_{2}a\right|^2 \right)
=\frac{1}{2} \left| \nabla a \right|^2
\text{,}
\end{equation*}
and we obtain 
\begin{equation*}
\sum\limits_{i=1}^{2} \left\langle \partial_{i}a, b\right\rangle^{2}
=|b|^2 \sum\limits_{i=1}^{2}
\left\langle \partial_{i}a, \frac{b}{|b|}\right\rangle^{2}
\leq
\frac{1}{2}\left| \nabla a\right|^{2}\left| b\right|^{2}
\text{.}
\end{equation*}
Therefore,
\begin{equation}
\begin{array}{l}
\underset{\mathbb{U}_r (w)}{\int }
\lambda \left( \rho- d^{\ast}_{k}(x) \right)
\cdot \left| \nabla G (x) \right| ^{2}
\text{ } dx
\\
=
\underset{\mathbb{U}_r (w)\cap \Omega^{\ast}_{k}(\rho)}{\int }
\lambda \left( \rho- |G^{\ast}_{k}-G(x)| \right)
\cdot \left| \nabla G (x) \right| ^{2}
\text{ } dx
\\
=
\underset{\mathbb{U}_r (w)\cap \Omega^{\ast}_{k}(\rho)}{\int }
\frac{\lambda^{\prime} \left( \rho- |G^{\ast}_{k}-G(x)| \right)}
{|G^{\ast}_{k}-G(x)|}
\cdot
\sum\limits_{i=1}^{2}
\left\langle \partial_{i} G (x), G^{\ast}_{k}-G(x) \right\rangle^{2}
\text{ } dx
\\
\leq
\frac{\rho}{2} \cdot
\underset{\mathbb{U}_r (w)\cap \Omega^{\ast}_{k}(\rho)}{\int }
\lambda^{\prime} \left( \rho- |G^{\ast}_{k}-G(x)| \right)
\cdot
\left| \nabla G (x) \right| ^{2}
\text{ } dx
\\
=
\frac{\rho}{2} \cdot
\underset{\mathbb{U}_r (w)}{\int }
\lambda^{\prime} \left( \rho- d^{\ast}_{k}(x) \right)
\cdot
\left| \nabla G (x) \right| ^{2}
\text{ } dx 
\end{array}
\label{eq:diff_ineq_1}
\end{equation}
where the inequality comes from applying the weak conformality conditions of
$G$ (proved in Proposition \ref{prop:wc}). 

Let
\begin{equation*}
\Psi_{k} (\rho)
:=
\underset{\mathbb{U}_r (w)}{\int }
\lambda \left( \rho- d^{\ast}_{k}(x) \right)
\cdot \left| \nabla G (x) \right| ^{2}
\text{ } dx
\text{.}
\end{equation*}
Then (\ref{eq:diff_ineq_1}) gives
\begin{equation*}
\Psi_{k} (\rho) \leq \frac{\rho}{2} \frac{d}{d\rho}\Psi_{k} (\rho) 
\end{equation*}
for all
$\rho \in (\rho_{k}, \frac{2}{5}\sigma_{k})$
for $k \in \{0,1,..., k_{0}-1\}$ or
$\rho \in (\rho_{k_{0}}, \frac{2}{5}\min\{\tau_{\ast},\sigma_{k_{0}}\})$.
This inequality implies the nondecreasing property of
$\frac{\Psi_{k} (\rho )}{\rho ^{2}}$,
i.e., if $s\le t$, then
\begin{equation}
\frac{\Psi_{k} (s)}{s^2} \le \frac{\Psi_{k} (t)}{t^2} 
\label{eq:mono1}
\end{equation}
where either $s, t \in (\rho_{k}, \frac{2}{5}\sigma_{k})$
as $k\in\{0,...,k_{0}-1\}$ or
$s, t \in (\rho_{k_{0}},\frac{2}{5}\min\{\tau_{\ast},\sigma_{k_{0}}\})$
as $k=k_{0}$.

\textbf{Step 3: Extending the monotonicity formula ``globally".}

%One might need to summary more efficiently the reason
%why we need to adjust the range of the parameter $\rho$.
%For example, besides the nested condition, one also needs
%to include the condition of keeping compact supported in perturbation.
The estimates of $\tau_{\ast}$ by the Dirichlet integral of $f$ rely on
applying the monotonicity formula (\ref{eq:mono1})
on each admissible closed ball
$\mathbb{B}^{\mathbf{Q}}_{s}(q^{(k)})$,
for each fixed $k\in\{0,...,k_{0}\}$,
and
keeping the nested condition,
\begin{equation}
\mathbb{B}^{\mathbf{Q}}_{s}(q^{(k)})
\subset
\mathbb{B}^{\mathbf{Q}}_{t}(q^{(k+1)}) 
\label{eq:nested-condition-1}
\end{equation}
by proper choices of $s$ and $t$,
\begin{equation}
\left(\rho_{k}, \frac{2}{5}\sigma_{k}\right)
\owns s<t\in
\left(\rho_{k+1}, \frac{2}{5}\min\{\tau_{\ast},\sigma_{k+1}\}\right)
\label{eq:nested-condition-2}
\end{equation}
for all $k\in\{0,...,k_{0}-1\}$.
However, one can't apply the monotonicity formula
of $\frac{\Psi_{k} (\rho)}{\rho^2}$ on the whole interval of
$\left(\rho_{k}, \frac{2}{5}\min\{\tau_{\ast},\sigma_{k}\}\right)$.
This is because that
the parameter $\rho$ in
$\lambda(\rho-d^{\ast}_{k}(\cdot))$
doesn't exactly represent the radius of an admissible closed ball
$\mathbb{B}^{Q}_{\rho}(\cdot)$
(note that there is an ``error" term, the harmonic function $h$,
in the definition of $d^{\ast}_{k}(\cdot)$).
Moreover, one needs to change the center of admissible closed balls,
i.e., from $q^{(k)}$ to $q^{(k+1)}$,
and keep the nested condition (\ref{eq:nested-condition-1})
to establish a relation
(which is nearly an inequality)
between $\frac{\Psi_{k} (s)}{s^2}$ and $\frac{\Psi_{k} (t)}{t^2}$
for $s$ and $t$ fulfilling (\ref{eq:nested-condition-2}).
But the parameter $\rho$ in
$\lambda(\rho-d^{\ast}_{k}(\cdot))$
only (nearly) represents the distance between
$f(x)$ and $q^{(k)}$ instead of $f(x)$ and $q^{(0)}$.
Thus, in order to keep the condition,
$\Omega^{\ast}_{k}(s)\subset\subset\mathbb{U}_{r}(w)$,
one needs to choose the upper bound of $s$ in (\ref{eq:nested-condition-2}).

Let $\lambda$ be the cut-off function defined in (\ref{eq:lambda-def}) and
$\varepsilon$ in (\ref{eq:lambda-def}) is any number satisfying
\begin{equation}
0<\varepsilon<\frac{\min\{\sigma_{0}, \tau_{\ast}\}}{10}
\text{.}
\label{eq:Upper_Bdd-varepsilon}
\end{equation}
Recall, from (\ref{eq:Lower_Bdd-sigma_(k)/rho_(k)}), that
$10Q \cdot \rho_{k}\le\sigma_{k}$ for each fixed $k\in\{0,...,L\}$
and, from (\ref{eq:def-k_0}), that $k_{0}\in\{0,...,L\}$ is the unique one fulfilling
$10 Q\cdot\rho_{k_{0}}\lvertneqq\tau_{\ast}\le 10 Q\cdot\rho_{k_{0}+1}$.
Let $\rho$ satisfy
$\varepsilon <\rho< \frac{2}{5}\min\{\sigma_{0},\tau_{\ast}\}$.
Note that, from (\ref{eq:Omega^ast_(k)subsetsubsetU_r}),
$\Omega^{\ast}_{0}(\rho)\subset\subset\mathbb{U}_{r}(w)$,
if $\rho$ fulfills the condition $\rho< \min\{\sigma_{0},\tau_{\ast}\}$.
Observe that
\begin{equation}
\begin{array}{l}
\rho^{-2}\cdot \Psi_{0} (\rho)
=
\rho^{-2}
\underset{\mathbb{U}_r (w)\cap \Omega^{\ast}_{0}(\rho)}{\int}
\lambda \left( \rho- d^{\ast}_{0}(x) \right)
\cdot \left| \nabla G (x) \right|^{2} \text{ } dx
\\ \\ \geq
\rho^{-2}
\underset{\mathbb{U}_r (w)\cap \Omega^{\ast}_{0}(\rho-\varepsilon)}{\int}
\lambda \left( \rho- d^{\ast}_{0}(x) \right)
\cdot \left| \nabla G (x) \right|^{2} \text{ } dx
=
\rho^{-2}
\underset{\mathbb{U}_r (w)\cap \Omega^{\ast}_{0}(\rho-\varepsilon)}{\int}
\left| \nabla G (x) \right|^{2} \text{ } dx
\text{.}
\label{eq:mono2-2_0}
\end{array}
\end{equation}
By applying the monotonicity formula (\ref{eq:mono1}) to the L.H.S. of
(\ref{eq:mono2-2_0})
and letting $\varepsilon\rightarrow 0$ on the R.H.S. of (\ref{eq:mono2-2_0}),
we derive
\begin{equation}
\rho^{-2}
\underset{\mathbb{U}_r (w)\cap \Omega^{\ast}_{0}(\rho)}
{\int }
\left| \nabla G (x) \right|^{2} \text{ } dx
\le
t^{-2}
\underset{\mathbb{U}_r (w)\cap \Omega^{\ast}_{0}(t)}
{\int }
\left| \nabla G (x) \right|^{2} \text{ } dx 
\label{eq:mono3-2_0}
\end{equation}
for $0 < \rho< t < \frac{2}{5}\min\{\sigma_{0},\tau_{\ast}\}$.
Now, since $w^{\ast}\in A$ is also a ``good" point of $G$
(see Remark \ref{rk:good_point_G}),
we may apply Lemma \ref{lem:Gr2.5} to the L.H.S. of (\ref{eq:mono3-2_0})
and let $t \rightarrow \frac{2}{5} \min\{\sigma_{0},\tau_{\ast}\}$
on the R.H.S. of (\ref{eq:mono3-2_0}) to derive
\begin{equation}
2 \pi
\le
\left( \frac{2}{5} \min\{\sigma_{0},\tau_{\ast}\} \right)^{-2}
\underset{\mathbb{U}_r(w)
\cap\Omega^{\ast}_{0}(\frac{2}{5}\min\{\sigma_{0},\tau_{\ast}\})}
{\int}
\left| \nabla G (x) \right| ^{2} \ dx
\text{.}
\label{eq:Lower_Bdd-sigma_0_0}
\end{equation}

{\bf Case $1^\circ$ $k_{0}=0$.}

As $\tau_{\ast}\le\sigma_{0}$,
we derive from (\ref{eq:Lower_Bdd-sigma_0_0}),
\begin{equation}
\tau_{\ast}
\le \frac{5}{2}\cdot
\sqrt{ \frac{Dir(G; \mathbb{U}_{r}(w))}{2\pi} }
\text{.}
\label{eq:Upper_Bdd-tau-1_0}
\end{equation}
As $\sigma_{0}<\tau_{\ast}\le 10Q\cdot\rho_{1}$,
we apply (\ref{eq:Upper_Bdd-rho_(k+1)/sigma_(k)}) in
Proposition \ref{prop:Nested_Adm_Closed_Ball}
(i.e., $\rho_{1} < C_{0}(n,Q) \cdot \sigma_{0}$),
and (\ref{eq:Lower_Bdd-sigma_0_0}) to derive
\begin{equation}
\tau_{\ast}<25 Q\cdot C_{0}(n,Q)\cdot
\sqrt{\frac{Dir(G; \mathbb{U}_{r}(w))}{2\pi}}
\text{.}
\label{eq:Upper_Bdd-tau-2}
\end{equation}
Thus, from (\ref{eq:Upper_Bdd-tau-1_0})
and (\ref{eq:Upper_Bdd-tau-2}),
we conclude that in this case, 
\begin{equation}
\mathcal{G}\left(f(w^{\ast}), f_{\vert \partial \mathbb{U}_{r}(w)} \right)
<
\left(\frac{5}{2}+25 Q\cdot C_{0}(n,Q)\right) \cdot
\sqrt{\frac{Dir(G; \mathbb{U}_{r}(w))}{2\pi}}
\text{.}
\label{eq:Upper_Bdd-tau-1}
\end{equation}

{\bf Case $2^\circ$ $k_{0}\ge 1$.}

We also need to adjust the range of
the parameter $\rho$ in the definition of distance function
$d^{\ast}_{k}(\cdot)$ to take care of the ``error" term
coming from the harmonic function $h$
in $d^{\ast}_{k}(\cdot)$.
Recall from (\ref{eq:(6)})
the conditions of nested admissible closed balls that,
for any $k\in\{1,..., L\}$,
\begin{equation}
\mathbb{B}^{\mathbf{Q}}_{s}
(q^{(k-1)})\subset\mathbb{B}^{\mathbf{Q}}_{t}(q^{(k)})
\text{, if }
s\le \sigma_{k-1} \text{ and } t\ge\rho_{k}
\text{.}
\label{eq:nested-condition-3}
\end{equation}
Thus, from (\ref{eq:nested-condition-3}),
$\mathcal{G}(q^{(k)},f(x)) \le \rho_{k}$,
if $x\in\Omega^{\ast}_{k-1}(\sigma_{k-1})$.
Recall from the definition of $\Omega^{\ast}_{k}(\cdot)$
that
if $x\in\Omega^{\ast}_{k-1}(\sigma_{k-1})$,
then
$\mathcal{G}(q^{(k-1)},f(x)) \le \sigma_{k-1}$
and $|h(w^{\ast})-h(x)|\le \sigma_{k-1}$.
Therefore, we obtain 
\begin{equation*}
d^{\ast}_{k}(x)\le\mathcal{G}(q^{(k)},f(x))
+|h(w^{\ast})-h(x)| \le \rho_{k}+\sigma_{k-1}
\text{, } \forall\text{ } x\in\Omega^{\ast}_{k-1}(\sigma_{k-1})
\text{.}
\end{equation*}
In other words,
\begin{equation}
s \le \sigma_{k-1} < \rho_{k}+\sigma_{k-1}\le t
\text{ }
\Rightarrow
\text{ }
\Omega^{\ast}_{k-1}(s)
\subset
\Omega^{\ast}_{k}(t)
\label{eq:Omega_(k-1)_subset_Omega_(k)}
\end{equation}
for all $k\in\{1,..., k_{0}\}$.

From the definition of $k_{0}$ and the assumption
$k_{0}\ge 1$, 
we have
$\sigma_{0}<\rho_{1}< 10Q\cdot \rho_{1}
\lvertneqq \tau_{\ast}$.
Thus, from (\ref{eq:Upper_Bdd-varepsilon}), 
we obtain 
$\varepsilon\in(0,\frac{\sigma_{0}}{10})$.
Note, from (\ref{eq:Omega^ast_(k)subsetsubsetU_r}), we have
\begin{equation*}
\Omega^{\ast}_{k}\left(\frac{2}{5}\min\{\tau_{\ast},\sigma_{k}\}\right)
\subset\subset\mathbb{U}_{r}(w)
\text{, }\forall \text{ } k\in\{0,...,k_{0}\}
\text{.}
%\label{eq:Omega_k-Omega_(k_0)}
\end{equation*}
From Proposition \ref{prop:Nested_Adm_Closed_Ball}
on the sequences of
$\{\rho_{k}\}_{k=0}^{L}$ and $\{\sigma_{k}\}_{k=0}^{L}$,
we have the inequality,
\begin{equation}
\rho_{k}+\sigma_{k-1}+\frac{\sigma_{0}}{10}
< 3Q \rho_{k}
< 4Q \rho_{k}
\le \frac{2}{5}\min\{\tau_{\ast},\sigma_{k}\}
\label{eq:Omega_k-Omega_(k_0)-1}
\end{equation}
for all $k\in\{1,..., k_{0}\}$.
Therefore, for any $k\in\{1,..., k_{0}\}$, 
$\rho_{k}+\sigma_{k-1}+\varepsilon<\frac{2}{5}\min\{\tau_{\ast},\sigma_{k}\}$, 
if $\varepsilon\in(0, \frac{\sigma_{0}}{10})$.

Now for each $k\in\{1,..., k_{0}\}$,
we consider
$\rho\in (
\rho_{k}+\sigma_{k-1}+\varepsilon, \frac{2}{5}\min\{\tau_{\ast},\sigma_{k}\})$ below.
From
\begin{equation}
\begin{array}{l}
\rho^{-2}\cdot \Psi_{k} (\rho)
=
\rho^{-2}
\underset{\mathbb{U}_r (w)\cap \Omega^{\ast}_{k}(\rho)}{\int}
\lambda \left( \rho- d^{\ast}_{k}(x) \right)
\cdot \left| \nabla G (x) \right|^{2} \text{ } dx
\\ \\ \geq
\rho^{-2}
\underset{\mathbb{U}_r (w)\cap \Omega^{\ast}_{k}(\rho-\varepsilon)}{\int}
\lambda \left( \rho- d^{\ast}_{k}(x) \right)
\cdot \left| \nabla G (x) \right|^{2} \text{ } dx
=
\rho^{-2}
\underset{\mathbb{U}_r (w)\cap \Omega^{\ast}_{k}(\rho-\varepsilon)}{\int}
\left| \nabla G (x) \right|^{2} \text{ } dx 
\label{eq:mono2-2_2}
\end{array}
\end{equation}
we again apply the monotonicity formula (\ref{eq:mono1}) to
the L.H.S. of
(\ref{eq:mono2-2_2})
and letting $\varepsilon\rightarrow 0$
on the R.H.S. of (\ref{eq:mono2-2_2})
to derive
\begin{equation}
\rho^{-2}
\underset{\mathbb{U}_r (w)\cap \Omega^{\ast}_{k}(\rho)}
{\int }
\left| \nabla G (x) \right|^{2} \text{ } dx
\le
t^{-2}
\underset{\mathbb{U}_r (w)\cap \Omega^{\ast}_{k}(t)}
{\int }
\left| \nabla G (x) \right|^{2} \text{ } dx 
\label{eq:mono3-2_2}
\end{equation}
for
$\rho_{k}+\sigma_{k-1}+\frac{\sigma_{0}}{10} < \rho< t
< \frac{2}{5}\min\{\tau_{\ast},\sigma_{k}\}$.
Now we let let $\rho=3\rho_{k}$
on the L.H.S. of (\ref{eq:mono3-2_2})
and $t \rightarrow \frac{2}{5}\min\{\tau_{\ast},\sigma_{k}\}$
on the R.H.S. of (\ref{eq:mono3-2_2}),
we have
\begin{equation}
\begin{array}{l}
\left( 3\rho_{k} \right)^{-2}
\underset{\mathbb{U}_r (w)\cap \Omega^{\ast}_{k}(3 \rho_{k})}
{\int }
\left| \nabla G (x) \right|^{2}
\text{ } dx
\\ \\
\le
\left(\frac{2}{5}\min\{\tau_{\ast},\sigma_{k}\}\right)^{-2}
\underset{\mathbb{U}_r (w)\cap \Omega^{\ast}_{k}
\left(\frac{2}{5}\min\{\tau_{\ast},\sigma_{k}\}\right)}
{\int}
\left| \nabla G (x) \right|^{2}
\text{ } dx
\text{.}
\label{eq:mono4-2_2}
\end{array}
\end{equation}
Thus, for each $k\in\{1,..., k_{0}\}$,
we have
\begin{equation}
\begin{array}{l}
\left(\frac{2}{5}\min\{\tau_{\ast},\sigma_{k}\}\right)^{-2}
\underset{\mathbb{U}_r (w)\cap \Omega^{\ast}_{k}
(\frac{2}{5}\min\{\tau_{\ast},\sigma_{k}\})}{\int}
\left| \nabla G (x) \right|^{2} \text{ } dx
\\ \\
\ge
\left(3\rho_{k}\right)^{-2}
\underset{\mathbb{U}_r (w)\cap \Omega^{\ast}_{k}(3\rho_{k})}
{\int}
\left| \nabla G (x) \right|^{2} \text{ } dx
\\ \\
\ge
\left(3C_{0}(n,Q)\cdot\sigma_{k-1}\right)^{-2}
\underset{\mathbb{U}_r (w)\cap \Omega^{\ast}_{k-1}
(\frac{2}{5}\sigma_{k-1})}{\int}
\left| \nabla G (x) \right|^{2} \text{ } dx 
\label{eq:mono-gap-1_1_0}
\end{array}
\end{equation}
where the first inequality comes from applying
(\ref{eq:mono4-2_2}) and the second inequality comes from applying
(\ref{eq:Omega_(k-1)_subset_Omega_(k)}),
(\ref{eq:Upper_Bdd-rho_(k+1)/sigma_(k)}).
Then, by applying (\ref{eq:mono-gap-1_1_0}) inductively on
$k=k_{0},...,1$,
we have
\begin{equation}
\begin{array}[t]{l}
\left(\frac{2}{5}\min\{\tau_{\ast},\sigma_{k_{0}}\}\right)^{-2}
\underset{\mathbb{U}_r (w)\cap \Omega^{\ast}_{k_{0}}
(\frac{2}{5}\min\{\tau_{\ast},\sigma_{k_{0}}\})}{\int}
\left| \nabla G(x) \right|^{2} \text{ } dx
\\ \\
\ge
\left(\frac{2}{15 C_{0}(n,Q)}\right)^{2 k_{0}-2}
\cdot \left(\frac{2}{5}\sigma_{0}\right)^{-2}
\underset{\mathbb{U}_r (w)\cap \Omega^{\ast}_{0}(\frac{2}{5}\sigma_{0})}{\int}
\left| \nabla G(x) \right|^{2} \text{ } dx
\text{.}
\label{eq:mono-gap-1_1}
\end{array}
\end{equation}
Therefore, by applying
(\ref{eq:Lower_Bdd-sigma_0_0}) to the R.H.S. of (\ref{eq:mono-gap-1_1}),
we have
\begin{equation}
%\begin{array}[t]{l}
\min\{\tau_{\ast},\sigma_{k_{0}}\}
\le
\frac{5}{2}
\left(\frac{15C_{0}(n,Q)}{2}\right)^{k_{0}-1}\cdot
\sqrt{\frac{Dir\left(G; \mathbb{U}_r(w)\right)}{2\pi}}
\text{.}
\label{eq:mono-gap-1_2}
%\end{array}
\end{equation}

As $\tau_{\ast}\le\sigma_{k_{0}}$,
we derive from (\ref{eq:mono-gap-1_2}) that
\begin{equation}
\tau_{\ast}
\le
\frac{5}{2}
\left(\frac{15C_{0}(n,Q)}{2}\right)^{k_{0}-1}\cdot
\sqrt{\frac{Dir\left(G; \mathbb{U}_r(w)\right)}{2\pi}}
\text{.}
\label{eq:Upper_Bdd-tau-1_2}
\end{equation}
As $\tau_{\ast}>\sigma_{k_{0}}$,
(\ref{eq:mono-gap-1_2}) gives
\begin{equation}
%\begin{array}[t]{l}
\sigma_{k_{0}}
\le
\frac{5}{2}
\left(\frac{15C_{0}(n,Q)}{2}\right)^{k_{0}-1}\cdot
\sqrt{\frac{Dir\left(G; \mathbb{U}_r(w)\right)}{2\pi}}
\text{.}
\label{eq:mono-gap-1_3}
%\end{array}
\end{equation}
Since the choice of $k_{0}$ implies $\tau_{\ast}\le 10Q\cdot\rho_{k_{0}+1}$,
from (\ref{eq:Upper_Bdd-rho_(k+1)/sigma_(k)}) in
Proposition \ref{prop:Nested_Adm_Closed_Ball}
(i.e., $\rho_{k+1} < C_{0}(n,Q) \cdot \sigma_{k}$)
and (\ref{eq:mono-gap-1_3}),
we have
\begin{equation}
\tau_{\ast}
\le
25Q \cdot C_{0}(n,Q)\cdot
\left(\frac{15 C_{0}(n,Q)}{2}\right)^{k_{0}-1}\cdot
\sqrt{\frac{Dir\left(G; \mathbb{U}_r(w)\right)}{2\pi}}
\text{.}
\label{eq:Upper_Bdd-tau-1_4}
\end{equation}

Since $C_{0}(n,Q)>1$ and $k_{0}\le Q-1$
(see Proposition \ref{prop:Nested_Adm_Closed_Ball}),
from
(\ref{eq:Upper_Bdd-tau-1}), (\ref{eq:Upper_Bdd-tau-1_2})
and (\ref{eq:Upper_Bdd-tau-1_4}),
we conclude that
\begin{equation}
\mathcal{G}\left(f(w^{\ast}), f_{\vert \partial \mathbb{U}_{r}(w)} \right)
=:\tau_{\ast}
\le 
\sqrt{\frac{Dir(G; \mathbb{U}_{r}(w))}{2\pi\cdot\delta(n,Q)}}
\label{eq:Upper_Bdd-tau-2_2}
\end{equation}
for some constant $\delta(n,Q)>0$.

\endproof

\subsection{\protect\smallskip Proof of Theorem \ref{thm:main}} 
To prove interior continuity of $f$, we may assume without loss of generality that
$\Omega=\mathbb{U}_{R_0}(0)\subset\mathbb{R}^2$, an open ball of radius $R_0>0$ 
with center at the origin of $\mathbb{R}^2$.
Since
$f\in\mathcal{Y}_{2}(\mathbb{U}_{R_0}(0),\mathbf{Q}_{Q}(\mathbb{R}^n))$
means
$\boldsymbol\xi\circ f \in W^{1, 2}(\mathbb{U}_{R_0}(0),\mathbb{R}^{N})$,
by Courant-Lebesgue Lemma (see Lemma \ref{lem:CL}),
for any $\mathbb{U}_R(w)\subset\subset\mathbb{U}_{R_0}(0)$,
one may choose a proper slice of $f$ by
$\partial\mathbb{U}_r(w)$ for some
$r\in\left[\frac{R}{2}, R\right]$
such that
$\boldsymbol\xi\circ f$ is continuous on the compact set
$\partial \mathbb{U}_r(w)$
and the oscillation of $\boldsymbol\xi\circ f$ is bounded by
$C(n, Q)\cdot \sqrt{Dir\left(\boldsymbol\xi\circ f; \mathbb{U}_{R}(w)\right)}$.
Thus, by the bi-Lipschitz continuity of $\boldsymbol\xi$,
we conclude that the multiple-valued function $f$ is uniformly continuous on
$\partial \mathbb{U}_r(w)$
and 
\begin{equation*}
\underset{\partial \mathbb{U}_{r}(w)}{\text{osc}} f
\leq
C(n, Q, \text{Lip}(\boldsymbol\xi), \text{Lip}(\boldsymbol\xi^{-1}))\cdot
\sqrt{Dir\left(f; \mathbb{U}_{R}(w)\right)}
=: \alpha_{1}(R)
\text{.}
%\label{eq:alpha_1}
\end{equation*}

On the other hand, 
for any ``good" point $y\in A\cap\mathbb{U}_{r}(w)$,
where $A$ is a set of Lebesgue points of
$|\nabla (\boldsymbol\xi_{0}\circ f)|^2$,
we apply Lemma \ref{lem:prethm} to derive
\begin{equation}
%\begin{array}{l}
\underset{x\in\partial\mathbb{U}_r(w)}
{\inf}
\mathcal{G}(f(x), f(y))
\le
\sqrt{\frac{Dir\left(G; \mathbb{U}_{r}(w)\right)}{2\pi\cdot\delta(n,Q)}}
\le
\frac{
\sqrt{Dir\left(f; \mathbb{U}_{R}(w)\right)}
+\lVert\nabla h\rVert_{L^{2}(\mathbb{U}_{R}(w))}
}
{\sqrt{2\pi\cdot\delta(n,Q)}}
\text{.}
%\end{array}
\label{eq:bdd_by_Dir(G)}
\end{equation}
Note that
$\lVert\nabla h\rVert_{L^{2}(\mathbb{U}_{R}(w))}$
in (\ref{eq:bdd_by_Dir(G)}) depends on the choice of
$\boldsymbol\xi_{0}$ because one needs to choose a proper coordinates of
$\mathbb{R}^{n}$ for a given $y$.
Therefore, we need to have a control
of $\lVert\nabla h\rVert_{L^{2}(\mathbb{U}_{R}(w))}$ in
(\ref{eq:bdd_by_Dir(G)})
when $h$ is induced from a distinct Lipschitz correspondence
$\boldsymbol\xi_{0}$.

Let $\overset{\star}{\boldsymbol\xi}_{0}$
be a \emph{fixed} Lipschitz correspondence
and
$\overset{\star}{h}:\mathbb{U}_{R_0}(0)\rightarrow\mathbb{R}^{2}$
be the harmonic function induced from the Hopf differential of
$\overset{\star}{\boldsymbol\xi}_{0}\circ f$. 
Suppose
$h:\mathbb{U}_{R_0}(0)\rightarrow\mathbb{R}^{2}$
is a harmonic function
with respect to an arbitrarily chosen Lipschitz correspondence
$\boldsymbol\xi_{0}$.
We would like to estimate the oscillation of
$
\lVert\nabla h \rVert_{L^{2}(\mathbb{U}_{R}(w))}
-\lVert\nabla \overset{\star}{h} \rVert_{L^{2}(\mathbb{U}_{R}(w))}
$ below.
By a simple computation from (\ref{eq:h=psi+bar(z)}),
we have
\begin{equation*}
|\nabla h|^2=\frac{|\varphi|^2}{8}+2
\text{.}
\end{equation*}
From applying
(\ref{eq:Upper_Bdd(varphi-bar_varphi)}) and (\ref{eq:|varphi|<|nabla_f|^2}),
we have
\begin{equation*}
\begin{array}{l}
\lvert |\nabla h|^2-|\nabla \overset{\star}{h}|^2 \rvert
\le
\frac{1}{8}\cdot |\varphi-\overset{\star}{\varphi}| \cdot
|\varphi+\overset{\star}{\varphi}|
\le
\frac{1}{8}\cdot
|\varphi-\overset{\star}{\varphi}|
\cdot (|\overset{\star}{\varphi}-\varphi| + 2|\overset{\star}{\varphi}|)
\\ \\
\le
2\cdot C_{R_0}^{2}
+ 2\cdot C_{R_0}\cdot
|\nabla(\overset{\star}{\boldsymbol\xi}_{0}\circ f)|^{2} 
\end{array}
\end{equation*}
where 
$C_{R_0}:=\frac{Dir(f; \mathbb{U}_{R_0}(0))}{\pi R_0^2}$.
Thus,
\begin{equation}
\begin{array}{l}
\int\limits_{\mathbb{U}_{R}(w)} \text{ } |\nabla h|^2 \text{ } dx
\le
\int\limits_{\mathbb{U}_{R}(w)} \text{ } |\nabla\overset{\star}{h}|^2 \text{ } dx
+ \int\limits_{\mathbb{U}_{R}(w)} \text{ }
\lvert |\nabla h|^2-|\nabla\overset{\star}{h}|^2 \rvert \text{ } dx
\\ \\
\le
Dir(\overset{\star}{h}; \mathbb{U}_{R}(w))
+ 2 \pi C_{R_0}^2 R^2
+ 2 C_{R_0}\cdot Dir(f; \mathbb{U}_{R}(w))
=: \beta(R)
\text{,}
\end{array}
\label{eq:beta(R)}
\end{equation}
where $\beta(R)\rightarrow 0$ as $R\rightarrow 0$.
Notice that
$\beta(R)$ defined on the R. H. S. of (\ref{eq:beta(R)}) is independent of
the choice of Lipschitz correspondence $\boldsymbol\xi_{0}$,
although
$\int\limits_{\mathbb{U}_{R}(w)} \text{ } |\nabla h|^2 \text{ } dx$
on the L. H. S. of (\ref{eq:beta(R)}) does depend on $\boldsymbol\xi_{0}$.
Thus, from (\ref{eq:bdd_by_Dir(G)}),
we obtain 
\begin{equation*}
%\begin{array}{l}
\underset{x\in \partial \mathbb{U}_r(w)}
{\inf}
\mathcal{G}(f(x),f(y))
\le
\frac{
\sqrt{Dir\left(f; \mathbb{U}_{R}(w)\right)}
+\sqrt{\beta(R)}
}
{\sqrt{2\pi\cdot\delta(n,Q)}}
=: \alpha_{2}(R) 
%\end{array}
\end{equation*}
for any $y\in \mathbb{U}_{r}(w)\cap A$.
Therefore,
\begin{equation*}
\underset{\mathbb{U}_{r}(w)\cap A}{\text{osc}}\text{ }
f
\le
4\max\{ \alpha_1 (R), \alpha_2 (R) \}
\text{.}
%\label{eq:osc_G_0}
\end{equation*}

Note that $\mathcal{L}^{2}(\mathbb{U}_{r}(w)\setminus A)=0$
and $\alpha_1 (R)$ and $\alpha_2 (R)$ tend to zero as $R\rightarrow 0$.
This proves the continuity of $f$ at $w$,
and the proof of Theorem \ref{thm:main} is finished.

{\bf Acknowledgement.}
The author started the project on stationary-harmonic multiple-valued functions in \cite{Lin01}. 
He would like to thank his thesis advisor
Bob Hardt for his warm encouragement and many helpful discussions at Rice University.
During the preparation of this manuscript, the author 
also received partial support
from the research grant of the National Science Council of Taiwan 
(NSC-100-2918-I-003-009), 
National Center for Theoretical Sciences in Taiwan
and the 
Max-Planck-Institute for Mathematics in the Sciences in Leipzig. 
%Max-Planck-Institut f\"{u}r Mathematik in den Naturwissenschaften, Leipzig. 
The author would like to acknowledge 
Professor Dr. Luckhaus, Professor Dr. Otto, Professor Dr. Stevens 
for their hospitality; and Dr. Spadaro for sharing opinion on the subject of 
Almgren's multiple-valued functions 
during his visiting in Leipzig.

\end{document}